\documentclass[a4paper,11pt]{article}

\usepackage{amsmath}
\usepackage{amsfonts}
\usepackage{amssymb}
\usepackage{amsthm}
\usepackage{graphicx}
\usepackage{psfrag}
\usepackage{subfigure}
\usepackage[german,english]{babel}
\usepackage[utf8]{inputenc}	
\usepackage{xparse}
\usepackage{mathtools}
\usepackage{dsfont}
\usepackage{enumerate}
\usepackage{hyperref}
\usepackage[backend=bibtex,style=numeric,maxnames=5]{biblatex}

\newtheorem{satz}{Satz}[section]
\newtheorem{theorem}[satz]{Theorem}
\newtheorem{lemma}[satz]{Lemma}

\theoremstyle{definition}
\newtheorem{definition}[satz]{Definition}

\newtheorem{bemerkung}[satz]{Remark}
\newtheorem{annahme}[satz]{Assumptions}
{\begin{proof}[Proof]}
{\end{proof}}

\makeatletter                    
\def\blfootnote{\xdef\@thefnmark{}\@footnotetext}
\DeclareMathAlphabet\mathbfcal{OMS}{cmsy}{b}{n}

\DeclarePairedDelimiter{\norm}{\lVert}{\rVert}
\NewDocumentCommand{\normH}{ s O{} m }{%
  \IfBooleanTF{#1}{\norm*{#3}}{\norm[#2]{#3}}_{L^2}%
}

\NewDocumentCommand{\normV}{ s O{} m }{%
  \IfBooleanTF{#1}{\norm*{#3}}{\norm[#2]{#3}}_{H^1}%
}
\renewcommand{\div}{\mathrm{div}}
\newcommand{\intO}{\int_{\Omega}}

\newcommand{\intT}{\int_{0}^{T}}

\newcommand{\ra}{\to}	
\newcommand{\N}{\mathds{N}}

\newcommand{\R}{\mathds{R}} 										
 
\newcommand{\bv}{\mathbf{v}}

\newcommand{\na}{\nabla}
\newcommand{\pa}{\partial}    

\renewcommand{\d}{\mathrm{d}}

\bibliography{Analysis_of_a_Cahn_Hilliard_Brinkman_model_for_tumour_growth_with_chemotaxis.bib}

\begin{document}
\begin{titlepage}
	\title{Analysis of a Cahn-Hilliard-Brinkman model for tumour growth with chemotaxis}
	\author{  Matthias Ebenbeck\footnote{Fakult\"at f\"ur Mathematik,  
			Universit\"at Regensburg,
			93040 Regensburg,
			Germany, e-mail: {\sf matthias.ebenbeck@mathematik.uni-regensburg.de}}\qquad\qquad
		Harald Garcke\footnote{Fakult\"at f\"ur Mathematik,  
			Universit\"at Regensburg,
			93040 Regensburg,
			Germany, e-mail: {\sf harald.garcke@mathematik.uni-regensburg.de}}}
		\date{}
\end{titlepage}
\maketitle
\begin{abstract}
	Phase field models recently gained a lot of interest in the context of tumour growth models. Typically Darcy-type flow models are coupled to Cahn-Hilliard equations. However, often Stokes or Brinkman flows are more appropriate flow models. We introduce and mathematically analyse a new Cahn-Hilliard-Brinkman model for tumour growth allowing for chemotaxis. Outflow boundary conditions are considered in order not to influence tumour growth by artificial boundary conditions. Existence of global-in-time weak solutions is shown in a very general setting.
\end{abstract}

\noindent{\bf AMS-Classification:} 
35K35, 
35Q92, 
92C50, 
35D30, 
76D07 


\section{Introduction}
\numberwithin{equation}{section}
Tumour growth models within the framework of continuum mechanics have been successful in describing many phenomena relevant for medical applications, see for example \cite{BellomoLiMaini,ByrnePreziosi,CristiniLowengrub,Friedman2,GarckeLamSitkaStyles,Grennspan,OdenTinsleyHawkins,RooseChapmanMaini}. First models based on differential equations focused on biochemical driving factors inhibiting or promoting the growth of the tumour. In the last twenty years also mechanical effects have been included in continuum mechanics based PDE modelling. Some simple models rely on a one species theory and more complex models use a multiphase mixture theory. All are based on fundamental balance laws and they differ in different mixture assumptions, in different constitutive laws involving the mechanical stresses or in different ways to account for cell adhesion mechanisms.\newline 
In the simplest case an equation of Darcy-type relating the velocity to the pressure gradient is proposed, see e.g. \cite{ByrneChaplain,CristiniLowengrubNie,FriedmanReitich,Grennspan}. This is motivated by the heterogeneous internal microstructure  of the tumour and in particular by the fact that the extracellular matrix (ECM) is considered as a porous media. Some models take aspects such as residual stress, plastic effects and vicoelasticity into account, see e.g. \cite{AmbrosiMollica,AmbrosiPreziosi}. Other authors, see e.g. \cite{FranksKing,Friedman,Friedman3}, use Stokes flow which can be motivated as an approximation of viscoelastic behaviour on sufficiently large time-scales. This seems to be a valid assumption as the time-scale of tumour growth is much larger than the relaxation times of viscoelastic biological materials, see \cite{AmbrosiPreziosi,FranksKing} for details. Cell-cell adhesion is typically modelled either within a sharp interface context leading to a free boundary problem involving the mean curvature of the interface, see \cite{Friedman,Friedman2, ZhengWiseChristini}, or within the context of phase field models of Cahn-Hilliard type, see \cite{ChristiniLiLowengrubWise,GarckeLamSitkaStyles,HilhorstKampmannNguyenZee}.\newline 
In the following, we will consider a Cahn-Hilliard-Brinkman system for tumour growth. For a bounded domain $\Omega\subset \R^d,~d=2,3,$ and a fixed time $T>0$, we consider for $Q\coloneqq \Omega\times (0,T)$ the following system of equations
\begin{subequations}
	\label{model_equations}
	\begin{alignat}{3}
	\label{model_eq_1}\div(\bv)&=\Gamma_{\textbf{v}}(\varphi,\sigma)&& \quad \text{in }   Q,\\
	\label{model_eq_2}-\div(T(\bv,p)) +\nu\bv&= \mu\na\varphi + (\chi_{\sigma}\sigma +\chi_{\varphi}(1-\varphi))\na\sigma && \quad \text{in } Q,\\
	\label{model_eq_3}\pa_t\varphi + \div(\varphi\bv) &= \div(m(\varphi)\na\mu)+\Gamma_{\varphi}(\varphi,\sigma,\mu)&& \quad \text{in } Q,\\
	\label{model_eq_4}\mu&= \epsilon^{-1}\psi'(\varphi)-\epsilon\Delta\varphi -\chi_{\varphi}\sigma && \quad \text{in } Q,\\
	\label{model_eq_5}\pa_t\sigma + \div(\sigma\bv) &= \div(n(\varphi)(\chi_{\sigma}\na\sigma- \chi_{\varphi}\na\varphi)) - \Gamma_{\sigma}(\varphi,\sigma,\mu)&& \quad \text{in } Q,
	\end{alignat}
\end{subequations}
where the viscous stress tensor is defined by
\begin{equation}
\label{definition_stress_tensor}T(\bv,p) = 2\eta(\varphi) D\bv+\lambda(\varphi)\div(\bv)\mathbf{I} - p\mathbf{I},
\end{equation}
and the symmetrised velocity gradient is given by
\begin{equation*}
D\bv\coloneqq \frac{1}{2}(\na\bv+\na\bv^T).
\end{equation*}
In (\ref{model_equations})-(\ref{definition_stress_tensor}), $\bv$ denotes the volume-averaged velocity of the mixture, $p$ denotes the pressure, $\sigma$ denotes the concentration of an unknown species acting as a nutrient, $\varphi\in[-1,1]$ denotes the difference in volume fractions,
with $\{\varphi = 1\}$ representing the unmixed tumour tissue, and  $\{\varphi = -1\}$ representing the surrounding healthy tissue, and $\mu$ denotes the chemical potential for $\varphi$. The functions $m(\cdot)$ and $n(\cdot)$ are positive functions representing the mobilities for the phase variable $\varphi$ and the nutrient density $\sigma$. The constant $\epsilon>0$ is related to the thickness of the diffuse interface. Moreover, the functions $\eta(\cdot)$ and $\lambda(\cdot)$ are non-negative and represent the shear and the bulk viscosity, respectively. The constants $\chi_{\sigma}$ and $\chi_{\varphi}$ are non-negative and related to the nutrient diffusion coefficient and the chemotaxis parameter. The right hand sides $\Gamma_{\bv}$, $\Gamma_{\varphi}$ and $\Gamma_{\sigma}$ account for volume changes due to growth, growth of the tumour, and sources respectively sinks for the nutrients.\newline
Furthermore, we define the free energy density of the nutrient by
\begin{equation}
\label{free_energy_kutrient}N(\varphi,\sigma) = \frac{\chi_{\sigma}}{2}|\sigma|^2 + \chi_{\varphi}\sigma(1-\varphi),
\end{equation}
and we denote the derivatives of $N$ by
\begin{equation*}
N_{\sigma}\coloneqq \frac{\pa N}{\pa\sigma} = \chi_{\sigma}\sigma + \chi_{\varphi}(1-\varphi),\quad N_{\varphi}\coloneqq \frac{\pa N}{\pa\varphi} = -\chi_{\varphi}\sigma.
\end{equation*}
Thus, we can rewrite (\ref{model_equations}) as 
\begin{alignat*}{3}
\div(\bv)&=\Gamma_{\textbf{v}}(\varphi,\sigma)&& \quad \text{in }   Q,\\
-\div(T(\bv,p)) +\nu\bv&= \mu\na\varphi + N_{\sigma}\na\sigma && \quad \text{in } Q,\\
\pa_t\varphi + \div(\varphi\bv) &= \div(m(\varphi)\na\mu)+\Gamma_{\varphi}(\varphi,\sigma,\mu)&& \quad \text{in } Q,\\
\mu&= \epsilon^{-1}\psi'(\varphi)-\epsilon\Delta\varphi +N_{\varphi} && \quad \text{in } Q,\\
\pa_t\sigma + \div(\sigma\bv) &= \div(n(\varphi)\na N_{\sigma}) - \Gamma_{\sigma}(\varphi,\sigma,\mu)&& \quad \text{in } Q.
\end{alignat*}
By $\mathbf{n}$ we will denote the outer unit normal on $\pa\Omega$, and $\pa_{\mathbf{n}} g \coloneqq \na g\cdot \mathbf{n}$ is the directional derivative.
We equip the system with the following initial and boundary conditions
\begin{subequations}
\label{boundary_and_initial_conditions}
\begin{alignat}{3}
\label{boundary_cond_1}\pa_{\mathbf{n}}\mu=\pa_{\mathbf{n}}\varphi &= 0 &&\quad \text{on }\pa\Omega\times (0,T)\eqqcolon\Sigma,\\
\label{boundary_cond_2}n(\varphi)\chi_{\sigma}\pa_{\mathbf{n}}\sigma &= b(\sigma_{\infty}-\sigma)&&\quad\text{on }\Sigma,\\
\label{boundary_cond_3}T(\bv,p)\mathbf{n} &= \mathbf{0}&&\quad\text{on }\Sigma,\\
\label{initial_cond}\varphi(0) = \varphi_0,\quad\sigma(0) &= \sigma_0 &&\quad\text{in }\Omega,
\end{alignat} 
\end{subequations}
where $\varphi_0,~\sigma_0,~\sigma_{\infty}$ are given functions and $b$ is a positive constant.

\subsection{Modelling aspects and comparison with other models}

In the following, we interpret and motivate the model under consideration. Furthermore, we give a comparison with other previous diffuse interface models in the literature.

\begin{enumerate}[$\bullet$]

\item The phase field variable $\varphi$ satisfies a convective Cahn-Hilliard-type equation with an additional source term $\Gamma_{\varphi}$, whereas the nutrient concentration is governed by a convection-diffusion-reaction equation. The fluxes in (\ref{model_eq_3}) and (\ref{model_eq_5}) are given by
\begin{align*}
\mathcal{J}_{\varphi}&\coloneqq -m(\varphi)\na\mu = -m(\varphi)\na(\epsilon^{-1}\psi'(\varphi)-\epsilon\Delta\varphi -\chi_{\varphi}\sigma),\\
\mathcal{J}_{\sigma}&\coloneqq -n(\varphi)\na (\chi_{\sigma}\sigma-\chi_{\varphi}\varphi).
\end{align*}
There are two non-standard terms in the definition of the fluxes. On the one hand, we have the term $m(\varphi)\na(\chi_{\varphi}\sigma)$, representing chemotactic mechanisms which drive the tumour cells towards regions of high nutrient. On the other hand, the term $n(\varphi)\na(\chi_{\varphi}\varphi)$ produces active transport mechanisms, which means that the nutrient cells are moving actively towards the tumour cells. For a detailed explanation of these mechanisms, we refer to \cite{GarckeLamSitkaStyles}.

\item When considering diffuse interface models for two-phase flows, there are two approaches that are commonly used to define the velocity. One possibility is to use a barycentric/mass-averaged velocity, leading to rather complicated expressions for the mass balance equations. We refer to the work of Lowengrub and Truskinovsky, cf. \cite{LowengrubTruskinovsky}, where the authors generalised the so-called \grqq Model H\grqq{}, see \cite{HalperinHohenberg}. In our model, we use a volume-averaged velocity of the form 
\begin{equation*}
\bv = u_1\bv_1 + u_2\bv_2,
\end{equation*}
where $u_i,\bv_i,i=1,2,$ are the volume fractions and the velocities of the fluids $i$.
This leads to a more simple expression for the balance equation (\ref{model_eq_1}), involving only the source terms and the densities of the single fluids. In this context, we refer to \cite{AbelsGarckeGrun,Boyer,GarckeLamSitkaStyles}. 
\item The source terms in the divergence and phase field equations are strongly related to each other. Indeed, denoting by $\Gamma_i,~ i=1,2,$ the source terms of the single components, one obtains
\begin{equation*}
\Gamma_{\varphi} = \frac{\Gamma_2}{\bar{\rho_2}} - \frac{\Gamma_1}{\bar{\rho_1}}\quad\text{and}\quad\Gamma_{\bv}=\frac{\Gamma_2}{\bar{\rho_2}}  +\frac{\Gamma_1}{\bar{\rho_1}},
\end{equation*}
where $\bar{\rho_1}$ and $\bar{\rho_2}$ are the mass densities of the pure components.
In the specific case that there is no loss or gain of mass locally, one obtains $\Gamma_2 = -\Gamma_1\eqqcolon \Gamma$. It was deduced in the work \cite{GarckeLamSitkaStyles} that
\begin{equation*}
\Gamma_{\bv} = c\Gamma_{\varphi},
\end{equation*}
where the constant $c$ depends only on the pure densities of the tumour and healthy components, $\bar{\rho_1}$ and $\bar{\rho_2}$. 
\item There are two choices for the source terms $\Gamma_{\varphi}$ and $\Gamma_{\sigma}$ which are commonly used in the literature. One possibility is to take
\begin{equation*}
\Gamma_{\varphi} = (P\sigma-A)h(\varphi),\quad \Gamma_{\sigma} = Ch(\varphi),
\end{equation*}
where $P$ is the proliferation rate, $A$ the apoptosis rate and $C$ is the consumption or nutrient uptake rate. Furthermore, $h$ is an interpolation function satisfying $h(-1)=0$ and $h(1)=1$ (e.g. $h(\varphi)=\min\{1,\max\{0,\frac{1}{2}(1+\varphi)\}\}$). These kind of source terms have been considered in \cite{ChristiniLiLowengrubWise,GarckeLam2,WiseLowengrubFrieboesChristini}. It is worth pointing out that in the tumour region $\{\varphi =1\}$, the growth of the tumour is proportional to the supply of nutrient, whereas in the healthy region $\{\varphi = -1\}$, the nutrient uptake is neglected due to the fact that the uptake of nutrient in the tumour region is much larger. \newline
Another possible choice for the source terms is given by
\begin{equation*}
\Gamma_{\varphi} = P(\varphi)(\sigma -\chi_{\varphi}\varphi-\mu),\quad \Gamma_{\sigma} = -\Gamma_{\varphi},
\end{equation*} 
with a non-negative proliferation function $P(\cdot)$. 
In \cite{HawkinsZeeKristofferOdenTinsley}, it has been suggested to take $P(\varphi) = p_0(1+\varphi)_{\text{+}}$. One could also choose $P(\varphi)= p_0(1-\varphi^2)\chi_{[-1,1]}(\varphi)$ with $\chi_{[-1,1]}$ being the characteristic function on the interval $[-1,1]$, as suggested in \cite{FrigeriGrasselliRocca} in the case $\chi_{\varphi}=1$. We also refer to \cite{ColliGilardiHilhorst} for the case $\chi_{\varphi}=0$.
\item A very important feature of our model is that the source term $\Gamma_{\bv}$ may depend on $\varphi$ and $\sigma$. Although this condition is of high practical relevance due to the relation between $\Gamma_{\bv}$ and $\Gamma_{\varphi}$, many authors have worked with prescribed source terms $\Gamma_{\bv}$ not depending on variables of the diffuse interface model, see e.g. \cite{GarckeLam1, JiangWuZheng}. This is related to the fact that boundary conditions of the form
\begin{equation*}
\bv = \mathbf{0}\text{ on }\pa\Omega\quad \text{or}\quad \bv\cdot\mathbf{n}=0\text{ on }\pa\Omega,
\end{equation*} 
require a source term $\Gamma_{\bv}$ which fulfils the compatibility condition:
\begin{equation*}
\intO \Gamma_{\bv}\d x = \intO \div(\bv)\d x = \int_{\pa\Omega}\bv\cdot\mathbf{n}\d \mathcal{H}^{d-1} = 0.
\end{equation*}
Also in the case of inhomogeneous boundary conditions in the form given above, a compatibility condition has to be satisfied. In the case of a solution dependent source term, it is in general not possible to fulfil such a condition. In the literature, there are only a few contributions in this direction, see e.g. \cite{GarckeLam4}, where they consider a quasi-static nutrient equation. Nevertheless, we have to assume that $\Gamma_{\bv}$ is a bounded function. Otherwise, we would have to estimate triple products of the form
\begin{equation*}
\intO \Gamma_{\bv}\mu\varphi\d x,\quad \intO \Gamma_{\bv}N_{\sigma}\sigma\d x,
\end{equation*} 
without having any a-priori-estimates on the solutions. However, in practice this does not lead to restrictions as $\varphi$ and $\sigma$ take bounded values in applications.
\item The energy of our model is given by
\begin{equation}
\label{energy}E(\varphi,\na\varphi,\sigma) = \intO \left(\frac{1}{\epsilon}\psi(\varphi)+\frac{\epsilon}{2}|\na\varphi|^2 + N(\varphi,\sigma)\right)\d x,
\end{equation}
where the first two terms describe the classical Cahn-Hilliard free energy. 
The last term is given by (\ref{free_energy_kutrient}) and consists of two parts.\newline 
The first term in (\ref{free_energy_kutrient}) leads to increasing energy of the total system generated by the presence of the nutrient. The second term in (\ref{free_energy_kutrient}) accounts for interactions between tumour and nutrient species. Indeed, for physical relevant values $\varphi\in [-1,1]$ and $\sigma\in [0,1]$, we observe that the second term attains its minimum when both $\varphi = +1$ (tumour region) and $\sigma = 1$. This results in chemotaxis and active transport mechanisms, driving the tumour cells towards regions with high nutrient supply and vice versa. In particular, without interaction the nutrient will only be driven by diffusion. For a more detailed motivation of the energy, we refer to \cite{GarckeLamSitkaStyles, HawkinsZeeKristofferOdenTinsley}.

\item The term $T(\bv,p)\mathbf{n}$ characterises effects due to friction on the boundary. Therefore, (\ref{boundary_cond_2}) can be referred to as a \grqq No-friction\grqq{} condition and is quite useful in applications, see \cite[App. III, 4.4]{Glowinski}. This condition is very popular for finite element discretizations of the Navier-Stokes equation since it appears naturally in the variational formulation of (\ref{model_eq_2}). In numerical simulations, it can be used to implement boundary conditions in an unbounded domain, for example a channel of infinite length. In this context, we also want to refer to the so-called classical \grqq Do-nothing\grqq{} boundary condition
\begin{equation}
\label{boundary_condition_do_kothing}\frac{\pa\bv}{\pa\mathbf{n}}-p\mathbf{n}=0,
\end{equation}
see e.g. \cite{HeywoodRannacherTurek}. Although (\ref{boundary_cond_3}) is of higher physical relevance, both (\ref{boundary_cond_3}) and (\ref{boundary_condition_do_kothing}) are less reflective than a Dirichlet boundary condition and therefore more useful in numerical applications.

\item Equation (\ref{model_eq_2}) is a Stokes-like equation, also referred to as the Brinkman equation when $\lambda(\cdot) \equiv 0$, $\eta(\cdot)\equiv \eta$ for a constant $\eta>0$ and $\Gamma_{\bv}=0$. The Brinkman model, which is a modification of Darcy's law, was first proposed by H.C. Brinkman in \cite{Brinkman} to model phase separation of isothermal, incompressible binary fluids in a porous media. This model has been analysed by several authors, see e.g. \cite{BosciaContiGrasselli,NgamsaadJirapornWannapong}.\newline
 The term $\mu\na\varphi + N_{\sigma}\na\sigma$ acts as a force in the momentum equation. Furthermore, the more general form of the stress tensor can be verified by the theory for isotropic, linearly viscous fluids, see \cite{EckGarckeKnabner, GurtinFriedAnand}.  The shear viscosity $\eta(\cdot)$ characterises the resistance of the fluid to shear, whereas the bulk viscosity $\lambda(\cdot)$ models the response of the fluid to changes in volume. 
 \newline
There are only a few contributions treating the case with variable viscosity, see  \cite{CollingsShenWise, ContiGiorgini, ContiGiorgini2, HuoZhangYang, ValdesParada}.
\end{enumerate}

\subsection{Notation and preliminaries}
We first want to fix some notation: For a (real) Banach space $X$ we denote by $\|.\|_X$ its norm, by $X^*$ the dual space and by $\langle .{,}. \rangle_X$ the duality pairing between $X^*$ and $X$. For an inner product space $X$, the inner product is denoted by $(.{,}.)_X$. We define the scalar product of two matrices by
\begin{equation*}
\mathbf{A}\colon \mathbf{B}\coloneqq \sum_{j,k=1}^{d}a_{jk}b_{jk}\quad\text{for } \mathbf{A},\mathbf{B}\in\R^{d\times d}.
\end{equation*}
 For the standard Lebesgue and Sobolev spaces with $1\leq p\leq \infty$, $k>0$,  we use the notation $L^p\coloneqq L^p(\Omega)$ and $W^{k,p}\coloneqq W^{k,p}(\Omega)$ with norms $\|.\|_{L^p}$ and $\|.\|_{W^{k,p}}$ respectively.  In the case $p=2$ we use $H^k\coloneqq W^{k,2}$ and the norm $\|.\|_{H^k}$. For $\beta\in (0,1)$ and $r\in (1,\infty)$, we will denote the Lebesgue and Sobolev spaces on the boundary by $L^p(\pa\Omega)$ and $W^{\beta,r}(\pa\Omega)$ with corresponding norms $\norm{\cdot}_{L^p(\pa\Omega)}$ and $\norm{\cdot}_{W^{\beta,r}(\pa\Omega)}$ (see \cite[Chap. I.3]{Sohr} for more details). By $\mathbf{L}^p$, $\mathbf{W}^{k,p}$, $\mathbf{H}^k$, $\mathbf{L}^p(\pa\Omega)$ and $\mathbf{W}^{\beta,r}(\pa\Omega)$, we will denote the corresponding spaces of vector valued and matrix valued functions. For the Bochner spaces, we use the notation $L^p(X)\coloneqq L^p(0,T;X)$ for a Banach space $X$ with $p\in [1,\infty]$.  For the dual space $X^*$ of a Banach space $X$, we introduce the (generalised) mean value by 
\begin{equation*}
v_{\Omega}\coloneqq \frac{1}{|\Omega|}\intO v\d x\quad\text{for } v\in L^1, \quad v_{\Omega}^*\coloneqq \frac{1}{|\Omega|}\langle v{,}1\rangle_X\quad\text{for } v\in X^*.
\end{equation*}
Moreover, we introduce the function spaces
\begin{align*}
&L_0^2\coloneqq \{w\in L^2\colon w_{\Omega}=0\},\quad H_N^2\coloneqq \{w\in H^2\colon \pa_{\mathbf{n}} w = 0 \text{ on } \pa\Omega\},\\
&(H^1)_0^*\coloneqq \{f\in (H^1)^*\colon f_{\Omega}^* =0\}.
\end{align*}
Then, the Neumann-Laplace operator $-\Delta_N\colon H^1\cap L_0^2\ra (H^1)_0^*$ is positive definite and self-adjoint. In particular, by the Lax-Milgram theorem and the Poincaré inequality (see (\ref{Poincare})), the inverse operator $(-\Delta_N)^{-1}\colon (H^1)_0^*\ra H^1\cap L_0^2$ is well-defined, and we set $u\coloneqq (-\Delta_N)^{-1}f$ for $f\in (H^1)_0^*$ if $u_{\Omega}=0$ and 
\begin{equation*}
-\Delta u = f\text{ in }\Omega,\quad\pa_{\mathbf{n}}u =0\text{ on }\pa\Omega.
\end{equation*}
We have dense and continuous embeddings $H_N^2\subset H^1\subset L^2\simeq (L^2)^*\subset (H^1)^*\subset (H_N^2)^*$ and the identifications $\langle u{,}v\rangle_{H^1}=(u{,}v)_{L^2}$, $\langle u{,}w\rangle _{H^2} = (u{,}w)_{L^2}$ for all $u\in L^2,~ v\in H^1$ and $w\in H_N^2$.\newline
We also want to recall Poincaré's inequality with mean value for $H^1$: There exists a constant $C_P$ depending only on $\Omega$ such that 
\begin{subequations}
\begin{equation}
\label{Poincare} \norm{f}_{L^2}\leq C_P(\norm{\na f}_{L^2} + |f_{\Omega}|)\quad\forall f\in H^1,
\end{equation}
or equivalently
\begin{equation}
\label{Poincare2}\norm{f-f_{\Omega}}_{L^2}\leq C_P\norm{\na f}_{L^2}\quad \forall f\in H^1.
\end{equation}
\end{subequations}
For convenience, we also recall Korn's inequality (see \cite[Thm. 6.3-3]{Ciarlet}): Let $\Omega\subset\R^d, d=2,3,$ be a bounded domain and $\mathbf{u}\in \mathbf{H}^1$. Then there exists a constant $C_K$ depending only on $\Omega$ such that
\begin{equation}
\label{Korn_inequality}\norm{\mathbf{u}}_{\mathbf{H}^1}\leq C_K\left(\norm{\mathbf{u}}_{\mathbf{L}^2}^2 + \intO D\mathbf{u}\colon D\mathbf{u}\d x\right)^{\frac{1}{2}}.
\end{equation}
We will also use the following Gronwall inequality in integral form, see \cite[Lemma 3.1]{GarckeLam3}: 

\begin{lemma}
	\label{Gronwall_lemma}
	Let $\alpha$, $\beta$, $u$ and $v$ be real-valued functions defined on $[0,T]$. Assume that $\alpha$ is integrable, $\beta$ is non-negative and continuous, $u$ is continuous, $v$ is non-negative and integrable. 
	If $u$ and $v$ satisfy the integral inequality
	\begin{equation*}
	u(s) + \int_{0}^{s}v(t)\d t\leq \alpha(s) + \int_{0}^{s}\beta(t)u(t)\d t\text{ for }s\in (0,T],
	\end{equation*}
	then it holds that for all $s\in (0,T]$
	\begin{equation}
	\label{Gronwall}u(s) + \int_{0}^{s}v(t)\d t\leq \alpha(s) + \int_{0}^{s}\alpha(t)\beta(t)\exp\!\bigg(\int_{0}^{t}\beta(r)\d r\bigg)\d t.
	\end{equation}	
\end{lemma}
Furthermore, we will use the following generalised Gagliardo-Nirenberg inequality:
\begin{lemma}\label{lemma_Gagliardo_Nirenberg}
Let $\Omega\subset \R^d,~d=2,3,$ be a bounded domain with Lipschitz boundary and $f\in W^{m,r}\cap L^q,~1\leq q,r\leq \infty$. For any integer $j,~0\leq j<m$, suppose there is $\alpha\in\R$ such that 
\begin{equation*}
j - \frac{d}{p} = \left(m-\frac{d}{r}\right)\alpha + (1-\alpha)\left(-\frac{d}{q}\right),\quad \frac{j}{m}\leq \alpha\leq 1.
\end{equation*}
Then, there exists a positive constant $C$ depending only on $\Omega,d,m,j,q,r,$ and $\alpha$ such that
\begin{equation}
\norm{D^jf}_{L^p}\leq C\norm{f}_{W^{m,r}}^{\alpha}\norm{f}_{L^q}^{1-\alpha}.
\end{equation}
\end{lemma}
The following interpolation inequality will also be of importance:
\begin{lemma}(\cite[Thm. II.4.1]{Galdi})\label{emma_interpolation_boundary}Let $\Omega\subset \R^d,d=2,3$ be a bounded domain with Lipschitz boundary and let $u\in W^{1,q}$ with $q\in [1,\infty)$. Assume
\begin{alignat*}{3}
r&\in [q,q(d-1)/(d-q)]&&\quad\text{if }q<d,\\
r&\in [q,\infty)&&\quad\text{if }q\geq d.
\end{alignat*}
Then, the following inequality holds:
\begin{equation}
\label{interpolation_inequality_boundary}\norm{u}_{L^r(\pa\Omega)}\leq C\left(\norm{u}_{L^q}^{1-\alpha}\norm{u}_{W^{1,q}}^{\alpha} + \norm{u}_{L^q}^{(1-\frac{1}{r})(1-\alpha)}\norm{u}_{W^{1,q}}^{\frac{1}{r}+\alpha(1-\frac{1}{r})}\right),
\end{equation}
where $C= C(d,r,q,\Omega)$ and $\alpha = d(r-q)/q(r-1)$.
\end{lemma}
We will also need the following theorem concerning solvability of the divergence equation:
\begin{lemma}( \cite[Sec. III.3]{Galdi})\label{lemma_divergence_equation}
Let $\Omega\subset\R^d,~d\geq 2,$ be a bounded domain with Lipschitz-boundary and let $1<q<\infty$. Then, for every $f\in L^q$ and $\mathbf{a}\in \mathbf{W}^{1-1/q,q}(\pa\Omega)$ satisfying
\begin{equation}
\label{divergence_compatibility_condition}\intO f\d x = \int_{\pa\Omega}\mathbf{a}\cdot\mathbf{n}\d \mathcal{H}^{d-1},
\end{equation}
there exists at least one solution $\mathbf{u}\in \mathbf{W}^{1,q}$ of the problem
\begin{alignat*}{3}
\div(\mathbf{u}) &= f&&\quad\text{in }\Omega,\\
\mathbf{u} &= \mathbf{a} &&\quad \text{on }\pa\Omega.
\end{alignat*}
In addition, the following estimate holds
\begin{equation}
\label{divergence_equation}\norm{\mathbf{u}}_{\mathbf{W}^{1,q}}\leq C(\norm{f}_{L^q}+\norm{\mathbf{a}}_{\mathbf{W}^{1-1/q,q}(\pa\Omega)}),
\end{equation}
with $C$ depending only on $\Omega$ and $q$.
\end{lemma}
Finally, in the Galerkin ansatz (see Sec. 3) we will make use of the following lemma (see \cite{AbelsTerasawa} for a proof):
\begin{lemma}\label{Stokes_Neumann result}
	Let $\Omega\subset\R^d, d=2,3,$ be a bounded domain with $C^{1,1}$-boundary and outer unit normal $\mathbf{n}$ and $1<q<\infty$. Furthermore, assume that $g\in W^{1,q}$, $\mathbf{f}\in \mathbf{L}^q$, $c\in W^{1,r}$ with $r>d$, and the functions $\eta(\cdot),~\lambda(\cdot)$ fulfil (A3) (see Assumptions \ref{Assumptions} below). Then, there exists a unique solution $(\bv,p)\in \mathbf{W}^{2,q}\times W^{1,q}$ of the system 
	\begin{subequations}
		\label{Stokes_subsystem}
		\begin{alignat}{3}\
		\label{Stokes_subsystem_1}-\div(2\eta(c) D\bv +\lambda(c)\div(\bv)\mathbf{I})+\nu\bv + \na p &= \mathbf{f}&&\quad \text{a.~e. in }\Omega,\\
		\label{Stokes_subsystem_2}\div(\bv) &= g&&\quad\text{a.~e. in }\Omega,\\
		\label{Stokes_subsystem_3}(2\eta(c) D\bv +\lambda(c)\div(\bv)\mathbf{I}-p\mathbf{I})\mathbf{n} &= \mathbf{0}&&\quad \text{a.~e. on }\pa \Omega,
		\end{alignat}
	\end{subequations}
	satisfying the following estimate
	\begin{equation}
	\label{Stokes_Neumann_estimate}\norm{\bv}_{\mathbf{W}^{2,q}} + \norm{p}_{W^{1,q}}\leq C(\norm{\mathbf{f}}_{\mathbf{L}^q}+\norm{g}_{W^{1,q}}),
	\end{equation}
	with a constant $C$ depending only on $\eta_0,~\eta_1,~\lambda_0,~q,~ \norm{c}_{W^{1,r}}$ and $\Omega$.
\end{lemma}
\begin{bemerkung}
	\begin{enumerate}
		\item[(i)] The statement also holds for a boundary of class $W^{2-\frac{1}{r},r}$ for some $r>d$. For a deeper discussion of this less restrictive condition, we refer to \cite{AbelsTerasawa}.
		\item[(ii)] Solutions of (\ref{Stokes_subsystem}) are stable under perturbations of $\mathbf{f},~g$ and $c$.
		\item[(iii)] In \cite{AbelsTerasawa}, they consider the case when $\lambda(\cdot) \equiv 0$ and with an inhomogeneous boundary condition in (\ref{Stokes_subsystem_3}). Using straightforward modifications, the result can be generalised to the case $\lambda(\cdot)\neq 0$. Indeed, using (\ref{Stokes_subsystem_2}), the terms involving $\lambda(\cdot)$ can be transformed into the r.h.s. of (\ref{Stokes_subsystem_1}) and (\ref{Stokes_subsystem_3}), respectively.
		\item[(iv)] Since $c\in W^{1,r}$ with $r>d$, the products $2\eta(c)D\bv$ and $\lambda(c)\div(\bv)\mathbf{I}$ belong to $\mathbf{W}^{1,q}$ for $\bv\in \mathbf{W}^{2,q}$. This results from the boundedness of the operator 
		\begin{equation*}
		\pi\colon W^{1,r}\times W^{1,q}\ra W^{1,q},\quad \pi(w,v)\mapsto wv,
		\end{equation*}
		for $1<q\leq r\leq \infty$, $r>d$, which is an easy consequence of the estimate
		\begin{equation*}
		\norm{wv}_{W^{1,q}}\leq C(\norm{w}_{L^{\infty}}\norm{v}_{W^{1,q}} + \norm{w}_{W^{1,r}}\norm{v}_{L^p}),\quad \frac{1}{q} = \frac{1}{p} + \frac{1}{r}
		\end{equation*}
		and the embeddings $W^{1,r}\hookrightarrow L^{\infty}$, $W^{1,r}\hookrightarrow W^{1,q}$, $W^{1,q}\hookrightarrow L^p$. The latter embedding follows since
		\begin{equation*}
		1-\frac{d}{q}\geq -\frac{d}{p}\Longleftrightarrow \frac{1}{d} - \frac{1}{q}\geq -\frac{1}{p} = \frac{1}{r} - \frac{1}{q}\Longleftrightarrow \frac{1}{d}\geq \frac{1}{r}\Longleftrightarrow r\geq d.
		\end{equation*}
	\end{enumerate}
\end{bemerkung}

\section{Main result}
We make the following assumptions:
\begin{annahme} \label{Assumptions}~
\begin{enumerate}
\item[(A1)] The constants $\epsilon, \chi_{\sigma}$ and $\nu$ are positive and fixed and $\chi_{\varphi}, b$ are fixed, non-negative constants.
\item[(A2)] The mobilities $m(\cdot),~n(\cdot)$ are continuous on $\R$ and satisfy
\begin{equation*}
m_0\leq m(t)\leq m_1,\quad n_0\leq n(t)\leq n_1\quad\forall t\in\R,
\end{equation*}
for positive constants $m_0,m_1,n_0,n_1$.
\item[(A3)] The viscosities fulfil $\eta,\lambda\in C^2(\R)$ with bounded first derivatives and
\begin{equation}
\label{assumption_on_viscosity}\eta_0\leq \eta(t)\leq \eta_1,\quad 0\leq \lambda(t)\leq \lambda_0\quad\forall t\in\R,
\end{equation}
for positive constants $\eta_0,\eta_1$ and a non-negative constant $\lambda_0$.
\item[(A4)] The functions $\Gamma_{\varphi}$ and $\Gamma_{\sigma}$ are of the form
\begin{align}
\nonumber\Gamma_{\varphi}(\varphi,\sigma,\mu) &= \Lambda_{\varphi}(\varphi,\sigma)-\theta_{\varphi}(\varphi,\sigma)\mu,\\
\label{assumption_source_terms_1}\Gamma_{\sigma}(\varphi,\sigma,\mu) &= \Lambda_{\sigma}(\varphi,\sigma)-\theta_{\sigma}(\varphi,\sigma)\mu,
\end{align}
where $\theta_{\varphi},\theta_{\sigma}\colon\R^2\ra\R$ are continuous bounded functions with $\theta_{\varphi}$ non-negative and \newline $\Lambda_{\varphi},\Lambda_{\sigma}\colon\R^2\ra\R$ are continuous with linear growth, i.~e.
\begin{equation}
\label{assumption_source_terms_2}|\theta_i(\varphi,\sigma)|\leq R_0,\quad |\Lambda_i(\varphi,\sigma)|\leq R_0(1+|\varphi|+|\sigma|)\quad\text{for }i\in\{\varphi,\sigma\},
\end{equation}
such that
\begin{equation}
\label{assumption_source_terms_3}|\Gamma_{\varphi}|+|\Gamma_{\sigma}|\leq R_0(1+|\varphi|+|\sigma|+|\mu|)
\end{equation}
for some positive constant $R_0$.
\item[(A5)] The function $\Gamma_{\textbf{v}}\in C^1(\R^2,\R)$ is assumed to be bounded, i.~e. 
\begin{equation}
\label{assumption_on_gamma_1}|\Gamma_{\textbf{v}}(\varphi,\sigma)|\leq \gamma_0,
\end{equation}
for a positive constant $\gamma_0$.
\item[(A6)] The function $\psi\in C^2(\R)$ is non-negative and satisfies 
\begin{equation}
\label{assumption_on_psi_1}\psi(t)\geq R_1|t|^2-R_2\quad\forall t\in\R
\end{equation}
for some positive constants $R_1,R_2,$ and either one of the following holds:
\begin{enumerate}[1.]
	\item If $\theta_{\varphi}$ is non-negative and bounded, then there exist positive constants $R_3,~R_4$ such that
	\begin{equation}
	\label{assumptions_on_psi_2}|\psi(t)|\leq R_3(1+|t|^2),\quad |\psi'(t)|\leq R_4(1+|t|),\quad |\psi''(t)|\leq R_4\quad\forall t\in\R.
	\end{equation}
	\item If $\theta_{\varphi}$ is positive and bounded, that is
	\begin{equation}
	\label{assumption_for_psi_3}R_0\geq \theta_{\varphi}(t,s)\geq R_5>0\quad\forall t,s\in\R,
	\end{equation}
	then
	\begin{equation}
	\label{assumption_on_psi_4}|\psi''(t)|\leq R_6(1+|t|^q),\quad q\in [0,4),
	\end{equation}
	for some positive constants $R_i,~i=5,6$. 
\end{enumerate}
Furthermore, we assume that
\begin{equation}
\label{assumption_on_epsilon}\frac{1}{\epsilon}> \frac{2\chi_{\varphi}^2}{\chi_{\sigma}R_1}.
\end{equation}
\item[(A7)] The initial and boundary data satisfy
\begin{equation}
\label{assumption_on_data}\sigma_{\infty}\in L^2(0,T;L^2(\pa\Omega)),\quad \varphi_0\in H^1,\quad\sigma_0\in L^2.
\end{equation}
\end{enumerate}
\end{annahme}
Due to the relation of $\epsilon$ to the thickness of the diffuse interface, which is typically very small, (\ref{assumption_on_epsilon}) in general means no restriction.\newline
We now introduce the weak formulation of (\ref{model_equations}), (\ref{boundary_and_initial_conditions}):
\begin{definition}(Weak solution for (\ref{model_equations}),~(\ref{boundary_and_initial_conditions})) \label{definition_weak_solution_1}We call a quintuple $(\varphi,\sigma,\mu,\bv,p)$ a weak solution of (\ref{model_equations}) and (\ref{boundary_and_initial_conditions}) if 
\begin{align*}
&\varphi \in L^{\infty}(0,T;H^1)\cap L^2(0,T;H^2)\cap W^{1,2}(0,T;(H^1)^*),\\
&\sigma\in L^{\infty}(0,T;L^2)\cap L^2(0,T;H^1)\cap W^{1,\frac{4}{3}}(0,T;(H^1)^*),\\
&\mu\in L^2(0,T;H^1),\quad \bv\in L^2(0,T;\mathbf{H}^1),\quad p\in L^{\frac{4}{3}}(0,T;L^2)
\end{align*}
such that 
\begin{align*}
&\div(\bv) = \Gamma_{\textbf{v}}\text{ a.~e. in }Q,\quad \varphi(0)=\varphi_0\text{ a.~e. in }\Omega,\\
&\langle \sigma(0){,}\zeta\rangle_{H^1} = \langle \sigma_0{,}\zeta\rangle_{H^1}\quad\forall \zeta\in H^1,
\end{align*}
and 
\begin{subequations}
\label{weak_form_1}
\begin{align}
\label{weak_form_1_eq_1}\intO T(\bv,p)\colon \na\boldsymbol{\Phi}+\nu\bv\cdot\boldsymbol{\Phi}\d x &= \intO (\mu\na\varphi + N_{\sigma}\na\sigma)\cdot\boldsymbol{\Phi}\d x,\\
\nonumber\langle\pa_t\varphi{,}\Phi\rangle_{H^1,(H^1)^*}  &= \intO -m(\varphi)\na\mu\cdot\na\Phi + \Gamma_{\varphi}\Phi\d x\\
\label{weak_form_1_eq_2} &\quad -\intO \na\varphi\cdot\bv \Phi  +\varphi\Gamma_{\textbf{v}}\Phi\d x,\\
\label{weak_form_1_eq_3}\intO \mu\Phi\d x &= \intO \epsilon^{-1}\Psi'(\varphi)\Phi + \epsilon\na\varphi\cdot\na\Phi - \chi_{\varphi}\sigma\Phi\d x,\\
\nonumber\langle\pa_t\sigma{,}\Phi\rangle_{H^1,(H^1)^*}  &= \intO -n(\varphi)\na N_{\sigma}\cdot\na\Phi - \Gamma_{\sigma}\Phi\d x\\
\label{weak_form_1_eq_4} &\quad -\intO \na\sigma\cdot\bv \Phi  +\sigma\Gamma_{\textbf{v}}\Phi\d x + \int_{\pa\Omega}b(\sigma_{\infty}-\sigma)\Phi,
\end{align}
\end{subequations}
for a.~e. $t\in(0,T)$ and for all $\boldsymbol{\Phi}\in \mathbf{H}^1,~\Phi\in H^1$.
\end{definition}

The main goal of this work is to prove the following existence result:
\begin{theorem}\label{main_theorem}(Weak solutions for (\ref{model_equations}), (\ref{boundary_and_initial_conditions}))
Let $\Omega\subset \R^d, d=2,3,$ be a bounded domain with $C^{1,1}$-boundary $\pa\Omega$. Suppose Assumption \ref{Assumptions} is satisfied. Then. there exists a weak solution quintuple $(\varphi,\sigma,\mu,\bv,p)$ for (\ref{model_equations}), (\ref{boundary_and_initial_conditions}) in the sense of Definition \ref{definition_weak_solution_1}.
Moreover, the following estimate holds:
\begin{align}
\nonumber &\norm{\varphi}_{L^{\infty}(H^1)\cap L^2(H^2)\cap W^{1,2}((H^1)^*)}+\norm{\sigma}_{L^{\infty}(L^2)\cap L^2(H^1)\cap W^{1,\frac{4}{3}}((H^1)^*)}\\
\nonumber &+\norm{\mu}_{L^2(H^1)}+b^{\frac{1}{2}}\norm{\sigma}_{L^2(L^2(\pa\Omega))} + \norm{p}_{L^{\frac{4}{3}}(L^2)}\\
\label{estimate_solution_1} &+\norm{\bv}_{L^2(\mathbf{H}^1)}+\norm{\div(\varphi\bv)}_{L^2(L^{\frac{3}{2}})}+\norm{\div(\sigma\bv)}_{L^{\frac{4}{3}}((H^1)^*)}\leq C,
\end{align}
for a constant $C$ depending only on the initial data, the domain $\Omega$ and the parameters of the system, but not on $(\varphi,\sigma,\mu,\bv,p)$. 
\end{theorem}

\section{Galerkin approximation} 
We will construct approximate solutions by applying a Galerkin approximation with respect to $\varphi,\mu$ and $\sigma$ and at the same time solve for $\bv$ and $p$ in the corresponding whole function spaces. As Galerkin basis for $\varphi,~\mu$ and $\sigma$, we will use the eigenfunctions of the Neumann-Laplace operator $\{w_i\}_{i\in\N}$ that form a basis of $L^2$. We will choose $w_1 = 1$.
By elliptic regularity, we see that $w_i\in H_N^2$ and for every $g\in H_N^2$ with $g_k\coloneqq\sum_{i=1}^{k}(g{,}w_i)_{L^2}w_i$ we obtain
\begin{equation*}
\Delta g_k = \sum_{i=1}^{k}(g{,}w_i)_{L^2}\Delta w_i = -\sum_{i=1}^{k}(g{,}\lambda_i w_i)_{L^2}w_i = \sum_{i=1}^{k}(g{,}\Delta w_i)_{L^2}w_i = \sum_{i=1}^{k}(\Delta g{,}w_i)_{L^2}w_i,
\end{equation*}
where $\lambda_i$ is the corresponding eigenvalue to $w_i$. Therefore, $\Delta g_k$ converges strongly to $\Delta g$ in $L^2$. Again using elliptic regularity theory, we obtain that $g_k$ converges strongly to $g$ in $H_N^2$. Thus the eigenfunctions $\{w_i\}_{i\in \N}$ of the Neumann-Laplace operator form an orthonormal Schauder basis in $L^2$ which is also a basis of $H_N^2$. \newline
We fix $k\in\N$ and define 
\begin{equation*}
\mathcal{W}_k\coloneqq \text{span}\{w_1,...,w_k\}.
\end{equation*}
Our aim is to find functions of the form 
\begin{equation*}
\varphi_k(t,x)=\sum_{i=1}^{k}a_i^k(t)w_i(x),~~\mu_k(t,x)=\sum_{i=1}^{k}b_i^k(t)w_i(x),~~\sigma_k(t,x)=\sum_{i=1}^{k}c_i^k(t)w_i(x)
\end{equation*}
satisfying the following approximation problem:
\begin{subequations}
	\begin{align}
	\label{approx_problem_eq_1a}\int_{\Omega}\pa_t\varphi_k v\d x &= \intO-\text{m}(\varphi_k)\na\mu_k\cdot\na v +\Gamma_{\varphi,k}v  -(\na\varphi_k\cdot \mathbf{v}_k+\varphi_k\Gamma_{\bv,k})v\d x,\\
	\label{approx_problem_eq_1b}\intO \mu_kv\d x &= \intO \epsilon\na\varphi_k\cdot \na v+ \epsilon^{-1}\psi'(\varphi_k)v - \chi_{\varphi} \sigma_kv\d x,\\
	\nonumber\intO\pa_t\sigma_kv\d x &= \intO -n(\varphi_k)(\chi_{\sigma}\na\sigma_k-\chi_{\varphi}\na\varphi_k)\cdot\na v-\Gamma_{\sigma,k}v - (\na\sigma_k\cdot\mathbf{v}_k+\sigma_k\Gamma_{\bv,k})v\d x\\
	\label{approx_problem_eq_1c}&\quad + \int_{\pa\Omega}b(\sigma_{\infty}-\sigma_k)v\d\mathcal{H}^{d-1},
	\end{align}
which has to hold for all $v\in\mathcal{W}_k$, where $\Gamma_{\varphi,k}\coloneqq \Gamma_{\varphi}(\varphi_k,\sigma_k,\mu_k),~\Gamma_{\sigma,k}\coloneqq \Gamma_{\sigma}(\varphi_k,\sigma_k,\mu_k)$ and $\Gamma_{\bv,k}\coloneqq \Gamma_{\bv}(\varphi_k,\sigma_k)$. Furthermore, we define the velocity $\bv_k$ and the pressure $p_k$ as the solutions of (\ref{Stokes_subsystem}) with 
	\begin{equation*}
	\mathbf{f} = \mu_k\na\varphi_k + N_{\sigma,k}\na\sigma_k,\quad g = \Gamma_{\bv,k}, \quad c = \varphi_k,
	\end{equation*}
	where
	\begin{equation*}
	N_{\sigma,k}\coloneqq \frac{\pa}{\pa \sigma}N(\varphi_k,\sigma_k) = \chi_{\sigma}\sigma_k + \chi_{\varphi}(1-\varphi_k).
	\end{equation*}
Using the continuous embedding $H_N^2\hookrightarrow L^{\infty}$ and (A5), straightforward arguments yield that
\begin{equation*}
\mu_k\na\varphi_k + N_{\sigma,k}\na\sigma_k\in L^2, \quad \Gamma_{\bv,k}\in H^1\cap L^{\infty}.
\end{equation*}
Therefore, by Lemma \ref{Stokes_Neumann result}, we obtain that $(\bv_k,p_k)\in \mathbf{H}^2\times H^1$ and the following equations are satisfied 
\begin{alignat}{3}
\label{approx_problem_eq_1d}-\div(T(\bv_k,p_k))+\nu\bv_k &= \mu_k\na\varphi_k + N_{\sigma,k}\na\sigma_k&&\quad\text{a.e. in }\Omega,\\
\label{approx_problem_eq_1e}\div(\bv_k) &= \Gamma_{\bv,k}&&\quad\text{a.e. in }\Omega,\\
\label{approx_problem_eq_1f}T(\bv_k,p_k)\mathbf{n} &= \mathbf{0}&&\quad\text{a.e. on }\pa\Omega.
\end{alignat}
\end{subequations}
We define the following matrices with components

\begin{equation*}
(\mathbf{S}_m^k)_{ji}\coloneqq \int_{\Omega}m(\varphi_k)\na w_i\cdot\na w_j\d x, \qquad (\mathbf{S}_n^k)_{ji}\coloneqq \int_{\Omega}n(\varphi_k)\na w_i\cdot\na w_j\d x~~\forall 1\leq i,j\leq n,
\end{equation*}
and introduce for all $1\leq i,j\leq n$ the notation

\begin{alignat*}{3}
\psi_j^k& \coloneqq \int_{\Omega}\psi'(\varphi_k)w_j\d x, \qquad&&\boldsymbol{\psi}^k \coloneqq (\psi_1^k,...,\psi_k^k)^T,\\
(\mathbf{M}_{\pa\Omega})_{ji}&\coloneqq \int_{\pa\Omega}w_iw_j\d \mathcal{H}^{d-1},\qquad&& \mathbf{S}_{ij}\coloneqq\int_{\Omega}\na w_i\cdot \na w_j\d x,\\
G_j^k&\coloneqq\int_{\Omega}\Gamma_{\varphi}(\varphi_k,\sigma_k,\mu_k)w_j,\qquad&& \mathbf{G}^k\coloneqq (G_1^k,...,G_k^k)^T,\\
F_j^k&\coloneqq\int_{\Omega}\Gamma_{\sigma}(\varphi_k,\sigma_k,\mu_k)w_j,\qquad&& \mathbf{F}^k\coloneqq (F_1^k,...,F_k^k)^T,\\
\Sigma_j^k&\coloneqq \int_{\pa\Omega}\sigma_{\infty}w_j\d \mathcal{H}^{d-1}\qquad &&\boldsymbol{\Sigma}^k\coloneqq (\Sigma_1^k,...,\Sigma_k^k)^T,\\
(\mathbf{C}^k)_{ji}&\coloneqq \intO \na w_i\cdot \textbf{v}_k w_j,\qquad &&(\mathbf{D}^k)_{ij}\coloneqq \intO w_i w_j\Gamma_{\bv}(\varphi_k,\sigma_k),
\end{alignat*}
and we denote by $\delta_{ij}$ the Kronecker-delta. Furthermore, we define $\mathbf{a}^k\coloneqq(a_1^k,...,a_k^k)^T$, $\mathbf{b}^k\coloneqq(b_1^k,...,b_k^k)^T$ and $\mathbf{c}^k\coloneqq (c_1^k,...,c_k^k)^T$. Inserting $v=w_j,~1\leq j\leq k,$ in (\ref{approx_problem_eq_1a})-(\ref{approx_problem_eq_1c}) and using the above introduced notation, we get a system of ODEs equivalent to (\ref{approx_problem_eq_1a})-(\ref{approx_problem_eq_1c}), given by

\begin{subequations}\label{ODE_system}
	\begin{align}
	\label{ODE_system_eq_1a}\frac{\d}{\d t}\mathbf{a}^k &= -\mathbf{S}_m^k\mathbf{b}^k + \mathbf{G}^k - (\mathbf{C}^k+\mathbf{D}^k)\mathbf{a}^k,\\
	\label{ODE_system_eq_1b}\mathbf{b}^k &= \epsilon \mathbf{S}\mathbf{a}^k + \epsilon^{-1}\boldsymbol{\psi}^k - \chi_{\varphi}\mathbf{c}^k,\\
	\label{ODE_system_eq_1c}\frac{\d}{\d t}\mathbf{c}^k &= \mathbf{S}_n^k(\chi_{\varphi}\mathbf{a}^k-\chi_{\sigma}\mathbf{c}^k) - \mathbf{F}^k - (\mathbf{C}^k+\mathbf{D}^k)\mathbf{c}^k - b\mathbf{M}_{\pa\Omega}\mathbf{c}^k + b\boldsymbol{\Sigma}^k ,
	\end{align} 
\end{subequations}
where $\bv_k,~p_k$ are defined as above. We complete the system with the following initial conditions:
\begin{subequations}\label{ODE_system_initial_conditions}
	\begin{align}
	\label{ODE_system_initial_cond_eq_1}(\mathbf{a}^k)_i(0)&=\int_{\Omega}\varphi_0w_i\d x~~\forall 1\leq i\leq k,\\
	\label{ODE_system_initial_cond_eq_2}(\mathbf{c}^k)_i(0)&=\int_{\Omega}\sigma_0w_i\d x~~\forall 1\leq i\leq k,
	\end{align}
\end{subequations}
where we have

\begin{align*}
\normV*{\sum_{i=1}^{k}(\mathbf{a}^k)_i(0)w_i}&\leq \normV{\varphi_0},\\
\normH*{\sum_{i=1}^{k}(\mathbf{c}^k)_i(0)w_i}&\leq \normH{\sigma_0}.
\end{align*}
Substituting (\ref{ODE_system_eq_1b}) and $\bv_k$ into (\ref{ODE_system_eq_1a}), (\ref{ODE_system_eq_1c}), we obtain a coupled system of ODEs for $\mathbf{a}^k$ and $\mathbf{c}^k$, where $\mathbf{S}_m^k,~\mathbf{S}_n^k,~\mathbf{C}^k$ and $\mathbf{D}^k$ depend non-linearly on the solutions $\mathbf{a}^k$ and $\mathbf{c}^k$. 
Owing to the continuity of $m(\cdot),~ n(\cdot),~ \psi'(\cdot),~ \Gamma_{\bv}(\cdot{,}\cdot)$ and the source terms and due to (A3) and the stability of the system (\ref{Stokes_subsystem}) under perturbations, we obtain that the r.h.s of (\ref{ODE_system}) depends continuously on $(\mathbf{a}^k,\mathbf{c}^k)$. \newline Therefore, the Cauchy-Peano theorem ensures that there exists $T_k^*\in (0,\infty]$ such that (\ref{ODE_system}), (\ref{ODE_system_initial_conditions}) has at least one solution triple $\mathbf{a}^k,\mathbf{b}^k,\mathbf{c}^k$ with $\mathbf{a}^k,\mathbf{b}^k,\mathbf{c}^k\in C^1([0,T_k^*),\R^k)$ (where we used the relation (\ref{ODE_system_eq_1b}) for $\mathbf{b}^k$).
Hence, (\ref{approx_problem_eq_1a})-(\ref{approx_problem_eq_1c}) admits at least one solution triplet $(\varphi_k,\mu_k,\sigma_k)\in C^1([0,T_k^*);\mathcal{W}_k)^3$.\newline Furthermore, we can define $\bv_k$ and $p_k$ as the solutions of (\ref{approx_problem_eq_1d})-(\ref{approx_problem_eq_1f}). With similar arguments as above, we obtain that $(\bv_k(t),p_k(t))\in \mathbf{H}^2\times H^1$ for all $t\in [0,T_k^*)$.

\section{A priori estimates}
In order to derive a-priori estimates, we will show energy estimates using the energy (\ref{energy}). However, the source terms $\Gamma_{\bv},~\Gamma_{\varphi}$ and $\Gamma_{\sigma}$ will make the a-priori estimates non-trivial.\newline
Let $\delta_{ij}$ denote the Kronecker-delta. We choose $v=b_j^kw_j$ in (\ref{approx_problem_eq_1a}), $v=\frac{\d}{\d t}a_j^kw_j$ in (\ref{approx_problem_eq_1b}) and $v=\chi_{\sigma}c_j^kw_j + \chi_{\varphi}(\delta_{1j}-a_j^k)w_j$ in (\ref{approx_problem_eq_1c}) and sum the resulting identities over $j=1,...,k,$ to obtain
\begin{align}
\nonumber &\frac{\d}{\d t} \intO \epsilon^{-1}\psi(\varphi_k)+\frac{\epsilon}{2}|\na\varphi_k|^2 + N(\varphi_k,\sigma_k)\d x\\
\nonumber &\quad  +\intO m(\varphi_k)|\na\mu_k|^2 + n(\varphi_k)|\na N_{\sigma,k}|^2\d x + \int_{\pa\Omega}b\chi_{\sigma}|\sigma_k|^2\d \mathcal{H}^{d-1}\\
\nonumber &\quad = \intO \Gamma_{\varphi,k}\mu_k -\Gamma_{\sigma,k}N_{\sigma,k}\d x+ \int_{\pa\Omega}b(\sigma_{\infty}N_{\sigma,k} - \sigma_k\chi_{\varphi}(1-\varphi_k))\d \mathcal{H}^{d-1}\\
\label{energy_identity_1} &\quad - \intO (\na\varphi_k\cdot\bv_k + \varphi_k \Gamma_{\bv,k})\mu_k + (\na\sigma_k\cdot\bv_k + \sigma_k\Gamma_{\bv,k})N_{\sigma,k}\d x.
\end{align}
where we used that
\begin{equation*}
N_{\varphi,k}\coloneqq \frac{\pa}{\pa\varphi}N(\varphi_k,\sigma_k) = -\chi_{\varphi}\sigma_k.   
\end{equation*}
For the Stokes subsystem, we would like to multiply (\ref{approx_problem_eq_1d}) with $\bv_k$ and integrate over $\Omega$. Then we would have to get an estimate for $p_k$ without having any a-priori-estimates on the solutions. Therefore, we use the so called method of subtracting the divergence. \newline 
Due to the assumptions on $\Omega$ and $\Gamma_{\bv}$ (in particular $\Gamma_{\bv,k}\in L^{\infty}$ for all $k\in\N$) and using Lemma \ref{lemma_divergence_equation}, there exists a solution $\mathbf{u}_k\in \mathbf{W}^{1,q},~q\in(1,\infty),$ (not necessarily unique) of the problem
\begin{alignat*}{3}
 \div(\mathbf{u}_k) &= \Gamma_{\bv,k}&& \quad\text{in }\Omega,\\
\mathbf{u}_k &= \frac{1}{|\pa\Omega|}\left(\intO\Gamma_{\bv,k}\d x\right)\mathbf{n}\eqqcolon \mathbf{a}_k&&\quad\text{on }\pa\Omega,
\end{alignat*}
 satisfying for every $q\in (1,\infty)$  the estimate
 \begin{equation}
 \label{divergence_equation_regularity}\norm{\mathbf{u}_k}_{W^{1,q}}\leq c\norm{\Gamma_{\bv,k}}_{L^q},
 \end{equation}
 with a constant $c$ depending only on $q$ and $\Omega$. We remark that (\ref{divergence_compatibility_condition}) is fulfilled since
 \begin{equation*}
\int_{\pa\Omega}\mathbf{a}_k\cdot\mathbf{n}\d \mathcal{H}^{d-1} = \frac{1}{|\pa\Omega|}\left(\intO\Gamma_{\bv,k}\d x\right)\int_{\pa\Omega}\mathbf{n}\cdot\mathbf{n}\d \mathcal{H}^{d-1} = \intO \Gamma_{\bv,k}\d x.
 \end{equation*}
Multiplying (\ref{approx_problem_eq_1d}) with $\bv_k-\mathbf{u}_k$, integrating over $\Omega$ and by parts and using (\ref{approx_problem_eq_1e})-(\ref{approx_problem_eq_1f}), we end up at
\begin{align}
\nonumber \intO 2\eta(\varphi_k)|D\bv_k|^2 +\nu|\bv_k|^2\d x &= \intO 2\eta(\varphi_k) D\bv_k\colon \na\mathbf{u}_k + \nu\bv_k\cdot \mathbf{u}_k\d x\\
\label{energy_equality_2}&\quad + \intO (\mu_k\na\varphi_k + N_{\sigma,k}\na\sigma_k)\cdot (\bv_k-\mathbf{u}_k).
\end{align}
Summing (\ref{energy_identity_1}) and (\ref{energy_equality_2}) gives
\begin{align}
\nonumber &\frac{\d}{\d t} \intO \epsilon^{-1}\psi(\varphi_k)+\frac{\epsilon}{2}|\na\varphi_k|^2 + N(\varphi_k,\sigma_k)\d x + \intO 2\eta(\varphi_k)|D\bv_k|^2 +\nu|\bv_k|^2\d x\\
\nonumber &\quad  +\intO m(\varphi_k)|\na\mu_k|^2 + n(\varphi_k)|\na N_{\sigma,k}|^2\d x + \int_{\pa\Omega}b\chi_{\sigma}|\sigma_k|^2\d \mathcal{H}^{d-1}\\
\nonumber &\quad = \intO \Gamma_{\varphi,k}\mu_k -\Gamma_{\sigma,k}N_{\sigma,k}\d x+ \int_{\pa\Omega}b(\sigma_{\infty}N_{\sigma,k} - \sigma_k\chi_{\varphi}(1-\varphi_k))\d \mathcal{H}^{d-1}\\
\nonumber &\quad - \intO (\na\varphi_k\cdot\mathbf{u}_k + \varphi_k \Gamma_{\bv,k})\mu_k + (\na\sigma_k\cdot\mathbf{u}_k + \sigma_k\Gamma_{\bv,k})N_{\sigma,k}\d x\\
 \label{energy_identity_3}&\quad +\intO 2\eta(\varphi_k) D\bv_k\colon \na\mathbf{u}_k + \nu\bv_k\cdot \mathbf{u}_k\d x.
\end{align}
\subsection{Estimation of the Stokes terms}
Using Hölder's and Young's inequalities, (\ref{assumption_on_viscosity}), (\ref{assumption_on_gamma_1}) and (\ref{divergence_equation_regularity}) with $q=2$, we see that
\begin{align}
\nonumber\left|\intO 2\eta(\varphi_k) D\mathbf{v}_k\colon\na\mathbf{u}_k + \nu\mathbf{v}_k\cdot\mathbf{u}_k\d x\right| &\leq 2\eta_1\norm{D\bv_k}_{\mathbf{L}^2}\norm{\na\mathbf{u}_k}_{\mathbf{L}^2} + \nu\norm{\bv_k}_{\mathbf{L}^2}\norm{\mathbf{u}_k}_{\mathbf{L}^2}\\
\label{A_priori_Stokes}&\leq \eta_0 \norm{D\bv_k}_{\mathbf{L}^2}^2 + \frac{\nu}{2}\norm{\bv_k}_{\mathbf{L}^2}^2 + C(q,|\Omega|)\left(\frac{\eta_1^2}{\eta_0}+\frac{\nu}{2}\right)\gamma_0^2,
\end{align}
where $C(q,|\Omega|)$ is the constant arising in (\ref{divergence_equation_regularity}).
\subsection{Estimation of the boundary term}
Using again Hölder's and Young's inequalities together with the trace theorem, we see that
\begin{align}
\nonumber &\left|\int_{\pa\Omega}b(\sigma_{\infty}N_{\sigma,k} - \sigma_k\chi_{\varphi}(1-\varphi_k))\d \mathcal{H}^{d-1}\right|\\
\nonumber &\leq \frac{b\chi_{\sigma}}{2}\norm{\sigma_k}_{L^2(\pa\Omega)}^2 + \left(\frac{2b\chi_{\varphi}^2}{\chi_{\sigma}}+\frac{b\chi_{\varphi}}{2}\right)(|\Omega|+\norm{\varphi_k}_{L^2(\pa\Omega)}^2) + b(\chi_{\varphi}+\chi_{\sigma})\norm{\sigma_{\infty}}_{L^2(\pa\Omega)}^2\\
\label{A_priori_boundary_integral}&\leq \frac{b\chi_{\sigma}}{2}\norm{\sigma_k}_{L^2(\pa\Omega)}^2 + C_1\left(1+\norm{\varphi_k}_{H^1}^2\right) + C_2\norm{\sigma_{\infty}}_{L^2(\pa\Omega)}^2,
\end{align}
where
\begin{equation*}
C_1\coloneqq \left(\frac{2b\chi_{\varphi}^2}{\chi_{\sigma}}+\frac{b\chi_{\varphi}}{2}\right) (|\Omega|+C_{\text{tr}}^2), \quad C_2\coloneqq b(\chi_{\varphi}+\chi_{\sigma}),
\end{equation*}
and $C_{\text{tr}}$ is the constant resulting from the trace theorem.

\subsection{Energy-inequality for non-negative $\theta_{\varphi}$}
First of all, we want to deduce an estimate for the $L^2$-norm of $\mu_k$. Inserting $v=b_j^kw_j$ into (\ref{approx_problem_eq_1b}) and summing over $j=1,...,k,$ yields
\begin{equation*}
\intO |\mu_k|^2\d x = \intO \epsilon^{-1}\psi'(\varphi_k)\mu_k + \epsilon\na\varphi_k\cdot\na\mu_k - \chi_{\varphi}\sigma_k\mu_k\d x.
\end{equation*} 
Using Hölder's and Young's inequalities together with the assumptions on $\psi$ (see (\ref{assumptions_on_psi_2})), we obtain
\begin{align*}
\norm{\mu_k}_{L^2}^2 &\leq \intO \epsilon^{-1}R_4(1+|\varphi_k|)|\mu_k| + \epsilon|\na\varphi_k||\na\mu_k| + \chi_{\varphi}|\sigma_k||\mu_k|\d x \\
&\leq \frac{1}{2}\norm{\mu_k}_{L^2}^2 + \frac{2R_4^2}{\epsilon^2}(|\Omega|+\norm{\varphi_k}_{L^2}^2) + \frac{\epsilon}{2}(\norm{\na\varphi_k}_{\mathbf{L}^2}^2+\norm{\na\mu_k}_{\mathbf{L}^2}^2) + \chi_{\varphi}^2\norm{\sigma_k}_{L^2}^2,
\end{align*}
and consequently
\begin{equation}
\label{A_priori_source_terms_b_1}\norm{\mu_k}_{L^2}^2\leq \frac{4R_4^2}{\epsilon^2}(|\Omega|+\norm{\varphi_k}_{L^2}^2) + \epsilon(\norm{\na\varphi_k}_{\mathbf{L}^2}^2+\norm{\na\mu_k}_{\mathbf{L}^2}^2) + 2\chi_{\varphi}^2\norm{\sigma_k}_{L^2}^2.
\end{equation}
By the specific form (\ref{assumption_source_terms_1}), we observe that
\begin{equation*}
\Gamma_{\varphi}(\varphi_k,\sigma_k,\mu_k)\mu_k = \Lambda_{\varphi}(\varphi_k,\sigma_k)\mu_k - \theta_{\varphi}(\varphi_k,\sigma_k)|\mu_k|^2.
\end{equation*}
Therefore, we can neglect the non-positive term $-\theta_{\varphi}(\varphi_k,\sigma_k)|\mu_k|^2$ on the r.h.s. of (\ref{energy_identity_3}). Using (\ref{assumption_source_terms_2}) and Hölder's inequality (in the following, we will write $\Lambda_{i,k}\coloneqq \Lambda_i(\varphi_k,\sigma_k)$ for $i=\varphi,\sigma$), we can estimate the first term on the r.h.s. of (\ref{energy_identity_3}) by
\begin{align*}
&\left|\intO \Lambda_{\varphi,k}\mu_k -\Gamma_{\sigma,k}(\chi_{\sigma}\sigma_k + \chi_{\varphi}(1-\varphi_k))\d x\right|\\
 &\quad \leq \normH{\Lambda_{\varphi,k}}\normH{\mu_k} + (\normH{\Lambda_{\sigma,n}}+ R_0\normH{\mu_k})(\normH{\chi_{\sigma}\sigma_k + \chi_{\varphi}(1-\varphi_k)})\\
 &\quad\leq R_0((1+\chi_{\varphi})(|\Omega|^{\frac{1}{2}}+\normH{\varphi_k})+(1+\chi_{\sigma})\normH{\sigma_k})\normH{\mu_k}\\
&\quad+ R_0(|\Omega|^{\frac{1}{2}}+\normH{\varphi_k}+\normH{\sigma_k})(\chi_{\varphi}|\Omega|^{\frac{1}{2}}+\chi_{\sigma}\normH{\sigma_k} + \chi_{\varphi}\normH{\varphi_k}).
\end{align*} 
Using Young's inequality, we obtain
\begin{equation}
\label{A_priori_source_terms_b_3}\left|\intO \Lambda_{\varphi,k}\mu_k -\Gamma_{\sigma,k}(\chi_{\sigma}\sigma_k + \chi_{\varphi}(1-\varphi_k))\d x\right| \leq \delta\normH{\mu_k}^2 +  C_{3,\delta}(1+\normH{\varphi_k}^2) + C_{4,\delta}\normH{\sigma_k}^2,
\end{equation}
with constants
\begin{align*}
C_{3,\delta}&\coloneqq \left(\frac{3R_0^2}{4\delta}(1+\chi_{\varphi})^2+R_0\left(1+\chi_{\varphi}+\chi_{\varphi}^2\right)\right)(1+|\Omega|),\\
C_{4,\delta}&\coloneqq \frac{3R_0^2}{4\delta}(1+\chi_{\sigma})^2+R_0 (1+\chi_{\sigma}+\chi_{\sigma}^2)
\end{align*}
and $\delta>0$ to be chosen later.
It remains to estimate the third and fourth integral on the r.h.s. of (\ref{energy_identity_3}). Using (\ref{assumption_on_gamma_1}), (\ref{divergence_equation_regularity}) and the continuous embedding $L^{\infty}\hookrightarrow L^q$ for all $q\in (1,\infty)$, we observe that
\begin{equation*}
\norm{\mathbf{u}_k}_{\mathbf{W}^{1,q}}\leq c(q,\Omega)\norm{\Gamma_{\bv,k}}_{L^q}\leq c(q,\Omega)\norm{\Gamma_{\bv,k}}_{L^{\infty}}\leq c(q,\Omega,\gamma_0),
\end{equation*}
for all $q\in (1,\infty)$.
Using Hölder's and Young's inequalities and the continuous embedding $W^{1,q}\hookrightarrow L^{\infty},~q\in (3,\infty)$, we obtain 
\begin{align}
\nonumber \left|\intO (\na\varphi_k\cdot\mathbf{u}_k +\varphi_k\Gamma_{\bv,k})\mu_k\d x\right| &\leq (\norm{\na\varphi_k}_{\mathbf{L}^2}\norm{\mathbf{u}_k}_{\mathbf{L}^{\infty}} + \normH{\varphi_k}\norm{\Gamma_{\bv,k}}_{L^{\infty}})\normH{\mu_k}\\
\nonumber & \leq C(q,|\Omega|)\norm{\Gamma_{\bv,k}}_{L^{\infty}}(\norm{\na\varphi_k}_{\mathbf{L}^2}+\normH{\varphi_k})\normH{\mu_k}\\
\label{A_priori_source_terms_b_5}& \leq \frac{\gamma_0^2 C(q,|\Omega|)}{2\delta}(\normH{\varphi_k}^2 + \norm{\na\varphi_k}_{\mathbf{L}^2}^2) + \delta\normH{\mu_k}^2,
\end{align}
for all $q\in (3,\infty)$ and with $\delta>0$ to be chosen later.
With similar arguments, we deduce for $q\in (3,\infty)$ that
\begin{equation}
\label{A_priori_source_terms_b_6}\left|\intO (\na\sigma_k\cdot\mathbf{u}_k + \sigma_k\Gamma_{\bv,k})N_{\sigma,k}\d x\right|\leq(C(q,|\Omega|)C_{5,\tilde{\delta}} + C_6)(1+\normH{\varphi_k}^2+\normH{\sigma_k}^2) + \tilde{\delta}\norm{\na\sigma_k}_{\mathbf{L}^2}^2,
\end{equation}
with 
\begin{equation*}
C_{5,\tilde{\delta}}\coloneqq \frac{\gamma_0^2(3\chi_{\varphi}^2(1+|\Omega|)+3\chi_{\sigma}^2)}{2\tilde{\delta}},\quad C_6\coloneqq \gamma_0\left(1+\chi_{\sigma} +\frac{\chi_{\varphi}}{2}(1+|\Omega|)\right),
\end{equation*}
and $\tilde{\delta}>0$ to be chosen later.
Furthermore, using Hölder's and Young's inequalities, we deduce that
\begin{equation}
\label{A_priori_source_terms_b_7}\norm{\chi_{\sigma}\na\sigma_k}_{\mathbf{L}^2}^2 = \norm{\na N_{\sigma,k}+\chi_{\varphi}\na\varphi_k}_{\mathbf{L}^2}^2
\leq 2(\norm{\na N_{\sigma,k}}_{\mathbf{L}^2}^2 + \norm{\chi_{\varphi}\na\varphi_k}_{\mathbf{L}^2}^2).
\end{equation}
In the following, we fix $q\in (3,\infty)$ and we denote by $C_K$ the constant arising in Korn's inequality. Choosing $\delta,\tilde{\delta}$ small enough and using (\ref{assumption_on_viscosity}), (\ref{assumption_on_psi_1}), (\ref{A_priori_Stokes})-(\ref{A_priori_source_terms_b_7}) in (\ref{energy_identity_3}), we obtain the following energy inequality
\begin{align}
\nonumber &\frac{\d}{\d t} \intO \epsilon^{-1}\psi(\varphi_k)+\frac{\epsilon}{2}|\na\varphi_k|^2 + \frac{\chi_{\sigma}}{2}|\sigma_k|^2 + \chi_{\varphi}\sigma_k(1-\varphi_k)\d x\\
\nonumber &\quad  + \frac{\min(\eta_0,\nu/2)}{C_K^2}\norm{\bv_k}_{\mathbf{H}^1}^2 + \frac{m_0}{2}\norm{\na\mu_k}_{\mathbf{L}^2}^2 + \frac{n_0\chi_{\sigma}^2}{2}\norm{\na\sigma_k}_{\mathbf{L}^2}^2 + \frac{b\chi_{\sigma}}{2}\norm{\sigma_k}_{L^2(\pa\Omega)}^2\\
\label{energy_identity_4}&\quad \leq \bar{C}_b(1+\norm{\na\varphi_k}_{\mathbf{L}^2}^2 + \normH{\sigma_k}^2 + \norm{\sigma_{\infty}}_{L^2(\pa\Omega)}^2) + \frac{\bar{C}_b}{R_1}(\norm{\psi(\varphi_k)}_{L^1} + R_2|\Omega|),
\end{align}
with a constant $\bar{C}_b$ depending on the system parameters, but not on $k\in\N$. 
Integrating with respect to time from $0$ to $s\in (0,T]$ gives
\begin{align}
\nonumber &\epsilon^{-1}\norm{\psi(\varphi_k(s))}_{L^1} + \frac{\epsilon}{2}\norm{\na\varphi_k(s)}_{\mathbf{L}^2}^2 + \frac{\chi_{\sigma}}{2}\normH{\sigma_k(s)}^2 + \intO \chi_{\varphi}\sigma_k(s)(1-\varphi_k(s))\d x\\
\nonumber &\quad + \int_{0}^{s}\frac{\min(\eta_0,\nu/2)}{C_K^2}\norm{\bv_k}_{\mathbf{H}^1}^2 + \frac{m_0}{2}\norm{\na\mu_k}_{\mathbf{L}^2}^2 + \frac{n_0\chi_{\sigma}^2}{2}\norm{\na\sigma_k}_{\mathbf{L}^2}^2 + \frac{b\chi_{\sigma}}{2}\norm{\sigma_k}_{L^2(\pa\Omega)}^2\d t\\
\nonumber &\quad\leq \bar{C}_b\left(1+\frac{R_2|\Omega|}{R_1}\right)T + \bar{C}_b\int_{0}^{s}\norm{\na\varphi_k}_{\mathbf{L}^2}^2 + \normH{\sigma_k}^2 + \frac{1}{R_1}\norm{\psi(\varphi_k)}_{L^1}\d t\\
\label{energy_identity_5}&\quad + \norm{\sigma_{\infty}}_{L^2(0,T;L^2(\pa\Omega))}^2 + \epsilon^{-1}\norm{\psi(\varphi_0)}_{L^1} + \frac{\epsilon}{2}\normV{\varphi_0}^2 + \frac{\chi_{\sigma}}{2}\normH{\sigma_0}^2.
\end{align}
Since $\varphi_0\in H^1,~\sigma_0\in L^2$ and $\psi(\varphi_0)\in L^1$ by assumption (\ref{assumptions_on_psi_2}), we observe that
\begin{equation*}
C_I\coloneqq \epsilon^{-1}\norm{\psi(\varphi_0)}_{L^1} + \frac{\epsilon}{2}\normV{\varphi_0}^2 + \frac{\chi_{\sigma}}{2}\normH{\sigma_0}^2< \infty.
\end{equation*}
Using Hölder's and Young's inequalities together with (\ref{assumption_on_psi_1}), (\ref{assumption_on_epsilon}), we obtain
\begin{align}
 \nonumber &\left|\intO \chi_{\varphi}\sigma_k(x,s)(1-\varphi_k(x,s))\d x\right| \\
 \nonumber&\quad \leq \frac{3\chi_{\sigma}}{8}\normH{\sigma_k(s)}^2 + \frac{\chi_{\varphi}^2}{\chi_{\sigma}R_1}\norm{\psi(\varphi_k(s))}_{L^1} + \left(\frac{\chi_{\varphi}^2R_2}{\chi_{\sigma}R_1} + \frac{2\chi_{\varphi}^2}{\chi_{\sigma}}\right)|\Omega|\\
 \label{A_priori_source_terms_b_8}&\quad\leq \frac{3\chi_{\sigma}}{8}\normH{\sigma_k(s)}^2 + \frac{1}{2\epsilon}\norm{\psi(\varphi_k(s))}_{L^1} + \left(\frac{\chi_{\varphi}^2R_2}{\chi_{\sigma}R_1} + \frac{2\chi_{\varphi}^2}{\chi_{\sigma}}\right)|\Omega|.
\end{align}
Substituting (\ref{A_priori_source_terms_b_8}) into (\ref{energy_identity_5}) yields
\begin{align}
\nonumber &\min\left\{\frac{1}{2\epsilon},\frac{\epsilon}{2},\frac{\chi_{\sigma}}{8}\right\}(\norm{\psi(\varphi_k(s))}_{L^1} + \norm{\na\varphi_k(s)}_{\mathbf{L}^2}^2 + \normH{\sigma_k(s)}^2) \\
\nonumber &\quad + \int_{0}^{s}\frac{\min(\eta_0,\nu/2)}{C_K^2}\norm{\bv_k}_{\mathbf{H}^1}^2 + \frac{m_0}{2}\norm{\na\mu_k}_{\mathbf{L}^2}^2 + \frac{n_0\chi_{\sigma}^2}{2}\norm{\na\sigma_k}_{\mathbf{L}^2}^2 + \frac{b\chi_{\sigma}}{2}\norm{\sigma_k}_{L^2(\pa\Omega)}^2\d t\\
\label{energy_identity_6} &\quad\leq \tilde{C}_b(1+T) + C_I+\int_{0}^{s}\tilde{C}_b(\norm{\na\varphi_k}_{\mathbf{L}^2}^2 + \normH{\sigma_k}^2 + \norm{\psi(\varphi_k)}_{L^1})\d t + \norm{\sigma_{\infty}}_{L^2(0,T;L^2(\pa\Omega))}^2,
\end{align}
where
\begin{equation*}
\tilde{C}_b\coloneqq \max\left\{\bar{C}_b\left(1+\frac{R_2|\Omega|}{R_1}\right), \left(\frac{\chi_{\varphi}^2R_2}{\chi_{\sigma}R_1} + \frac{2\chi_{\varphi}^2}{\chi_{\sigma}}\right)|\Omega|,\frac{\bar{C}_b}{R_1} \right\}.
\end{equation*}
Setting
\begin{equation*}
\alpha\coloneqq \tilde{C}_b(1+T) + C_I + \norm{\sigma_{\infty}}_{L^2(0,T;L^2(\pa\Omega))}^2,\quad\beta\coloneqq \tilde{C}_b, 
\end{equation*}
and noting that
\begin{equation*}
\alpha\left(1+\int_{0}^{s}\beta\exp\left(\int_{0}^{t}\beta\d r\right)\d t\right) = \alpha(1+\exp(\beta s)-1)\leq \alpha\exp(\beta T),
\end{equation*}
an application of Lemma \ref{Gronwall_lemma} to (\ref{energy_identity_6}) gives
\begin{align}
\nonumber&\sup_{s\in (0,T]}(\norm{\psi(\varphi_k(s))}_{L^1} + \norm{\na\varphi_k(s)}_{\mathbf{L}^2}^2 + \normH{\sigma_k(s)}^2)\\
\label{central_estimate_1}&\quad +\int_{0}^{T}\norm{\bv_k}_{\mathbf{H}^1}^2 + \norm{\na\mu_k}_{\mathbf{L}^2}^2 + \norm{\na\sigma_k}_{\mathbf{L}^2}^2 + \norm{\sigma_k}_{L^2(\pa\Omega)}^2\d t\leq C,
\end{align}
with a constant $C$ depending only on $T$ and the system parameters, but not on $k\in\N$. In particular, we remark that $C$ does not depend on $b$. In the following, we will use the constant $C$ (only depending on the system parameters and $T$, but not on $k\in\N$ and $b$) as a generic constant which may change even within one line. Using assumption (\ref{assumption_on_psi_1}) and (\ref{A_priori_source_terms_b_1}), an immediate consequence of (\ref{central_estimate_1}) is given by
\begin{equation}
\label{central_estimate_2}\sup_{s\in (0,T]}\normV{\varphi_k(s)} + \intT\normV{\mu_k}^2\leq C.
\end{equation}

\subsection{Energy-inequality for positive $\theta_{\varphi}$}
We assume that the assumptions (\ref{assumption_for_psi_3})-(\ref{assumption_on_psi_4}) for $\theta_{\varphi}$ and $\psi$ are valid. Then, arguing as above, the specific form of $\Gamma_{\varphi}$ yields
\begin{equation*}
\Gamma_{\varphi,k}\mu_k = \Lambda_{\varphi,k}\mu_k - \theta_{\varphi}(\varphi_k,\sigma_k)|\mu_k|^2.
\end{equation*}
We move the second term on the r.h.s. of this equation to the l.h.s. of (\ref{energy_identity_3}). Then, we can perform exactly the same estimates as in the last subsection, except from (\ref{A_priori_source_terms_b_1}). We remark that estimate (\ref{A_priori_source_terms_b_1}) was the only reason why we needed assumption (\ref{assumptions_on_psi_2}). Again chosing $\delta$ and $\tilde{\delta}$ small enough, we arrive at the following inequality (compare (\ref{energy_identity_6}))
\begin{align}
\nonumber &\min\left\{\frac{1}{2\epsilon},\frac{\epsilon}{2},\frac{\chi_{\sigma}}{8}\right\}(\norm{\psi(\varphi_k(s))}_{L^1} + \norm{\na\varphi_k(s)}_{\mathbf{L}^2}^2 + \normH{\sigma_k(s)}^2) \\
\nonumber &\quad + \int_{0}^{s}C_{11}\norm{\bv_k}_{\mathbf{H}^1}^2 + m_0\norm{\na\mu_k}_{\mathbf{L}^2}^2 + \frac{R_5}{2}\normH{\mu_k}^2 \frac{n_0\chi_{\sigma}^2}{2}\norm{\na\sigma_k}_{\mathbf{L}^2}^2 + \frac{b\chi_{\sigma}}{2}\norm{\sigma_k}_{L^2(\pa\Omega)}^2\d t\\
\label{energy_identity_6_positive} &\quad\leq C(1+T) + C_I+\int_{0}^{s}C(\norm{\na\varphi_k}_{\mathbf{L}^2}^2 + \normH{\sigma_k}^2 + \norm{\psi(\varphi_k)}_{L^1})\d t+ \norm{\sigma_{\infty}}_{L^2(0,T;L^2(\pa\Omega))}^2 ,
\end{align}
with $C_I$ as defined in the last subsection.
The reason why we have the term $m_0\norm{\na\mu_k}_{\mathbf{L}^2}^2$ instead of $\frac{m_0}{2}\norm{\na\mu_k}_{\mathbf{L}^2}^2$ is that we do not use (\ref{A_priori_source_terms_b_1}). Notice that we still have 
\begin{equation*}
\norm{\psi(\varphi_0)}_{L^1}<\infty,
\end{equation*}
since
\begin{equation*}
\norm{\psi(\varphi_0)}_{L^1}\leq C(1+\norm{\varphi_0}_{L^6}^6)\leq C(1+\normV{\varphi_0}^6)<\infty
\end{equation*}
due to assumption (\ref{assumption_on_psi_4}) and the continuous embedding $H^1\hookrightarrow L^6$.\newline
Again applying Lemma \ref{Gronwall_lemma}, from (\ref{energy_identity_6_positive}) we obtain
\begin{align}
\nonumber&\sup_{s\in (0,T]}(\norm{\psi(\varphi_k(s))}_{L^1} + \norm{\na\varphi_k(s)}_{\mathbf{L}^2}^2 + \normH{\sigma_k(s)}^2)\\
\label{central_estimate_3}&\quad +\int_{0}^{T}\norm{\bv_k}_{\mathbf{H}^1}^2 + \normV{\mu_k}^2 + \norm{\na\sigma_k}_{\mathbf{L}^2}^2 + \norm{\sigma_k}_{L^2(\pa\Omega)}^2\d t\leq C.
\end{align}

With similar arguments as above, assumption (\ref{assumption_on_psi_1}) yields
\begin{equation}
\label{central_estimate_4}\sup_{s\in (0,T]}\normV{\varphi_k(s)}\leq C.
\end{equation}

\subsection{Estimation of the pressure}
Taking the scalar product of (\ref{approx_problem_eq_1d}) with $\boldsymbol{\Phi}\in \mathbf{H}^1$, integrating over $\Omega$ and by parts, when using (\ref{approx_problem_eq_1d})-(\ref{approx_problem_eq_1f}) we obtain
\begin{align}
\nonumber\intO  p_k\div(\boldsymbol{\Phi})\d x &= \intO (2\eta(\varphi_k) D\bv_k + \lambda(\varphi_k)\Gamma_{\bv,k}\mathbf{I})\colon \na\boldsymbol{\Phi}\d x\\ 
\label{A_priori_pressure_1}&\quad+ \intO (\nu\bv_k - \mu_k\na\varphi_k - N_{\sigma,k}\na\sigma_k)\cdot\boldsymbol{\Phi}\d x 
\end{align}
for all $\boldsymbol{\Phi}\in \mathbf{H}^1$. Now, we define a family of functionals on $\mathbf{H}^1$ by
\begin{equation*}
\label{A_priori_pressure_2}\mathcal{F}_k(\boldsymbol{\Phi})\coloneqq \intO (2\eta(\varphi_k) D\bv_k + \lambda(\varphi_k)\Gamma_{\bv,k}\mathbf{I})\colon \na\boldsymbol{\Phi} + \nu\bv_k\cdot \boldsymbol{\Phi} - (\mu_k\na\varphi_k + N_{\sigma,k}\na\sigma_k)\cdot\boldsymbol{\Phi}\d x
\end{equation*}
for all $\boldsymbol{\Phi}\in \mathbf{H}^1$. Using Hölder's inequality, (\ref{assumption_on_viscosity}), (\ref{assumption_on_gamma_1}) and the continuous embedding $H^1\hookrightarrow L^6$, we obtain
\begin{equation*}
\nonumber |\mathcal{F}_k(\boldsymbol{\Phi})|\leq C(1+\norm{\bv_k}_{\mathbf{H}^1} + \norm{\mu_k}_{L^3}\norm{\na\varphi_k}_{\mathbf{L}^2} + \norm{N_{\sigma,k}}_{L^3}\norm{\na\sigma_k}_{\mathbf{L}^2})\norm{\boldsymbol{\Phi}}_{\mathbf{H}^1},
\end{equation*}
with $C = C(\Omega,\gamma_0,\eta_1,\lambda_0,\nu)$. Taking the supremum over all $\boldsymbol{\Phi}\in \mathbf{H}^1$ with $\norm{\boldsymbol{\Phi}}_{\mathbf{H}^1}\leq 1$, we deduce that
\begin{equation}
\label{A_priori_pressure_3} \norm{\mathcal{F}_k}_{(\mathbf{H}^1)^*}\leq  C(1+\norm{\bv_k}_{\mathbf{H}^1} + \norm{\mu_k}_{L^3}\norm{\na\varphi_k}_{\mathbf{L}^2} + \norm{N_{\sigma,k}}_{L^3}\norm{\na\sigma_k}_{\mathbf{L}^2}).
\end{equation}
From (\ref{A_priori_pressure_1}), we see that
\begin{equation}
\label{A_priori_pressure_4}\mathcal{F}_k(\boldsymbol{\Phi}) = \intO p_k\div(\boldsymbol{\Phi})\d x\quad\forall \boldsymbol{\Phi}\in \mathbf{H}^1.
\end{equation}
Now, using Lemma \ref{lemma_divergence_equation}, we deduce that there is at least one solution $\mathbf{q}_k\in \mathbf{H}^1$ of the system
\begin{alignat*}{3}
\div(\mathbf{q}_k) &= p_k&&\quad\text{in }\Omega,\\
\mathbf{q}_k &= \frac{1}{|\pa\Omega|}\left(\intO p_k\d x\right)\mathbf{n}&&\quad\text{on }\pa\Omega
\end{alignat*}
such that
\begin{equation}
\label{A_priori_pressure_5}\norm{\mathbf{q}_k}_{\mathbf{H}^1}\leq C_d\normH{p_k},
\end{equation}
with $C_d$ depending only on $\Omega$. Notice that the compatibility condition (\ref{divergence_compatibility_condition}) is satisfied since
\begin{equation*}
\int_{\pa\Omega} \mathbf{q}_k\cdot\mathbf{n}\d \mathcal{H}^{d-1} = \frac{1}{|\pa\Omega|}\left(\intO p_k\d x\right)\int_{\pa\Omega}\mathbf{n}\cdot\mathbf{n}\d \mathcal{H}^{d-1} = \intO p_k\d x.
\end{equation*}
Choosing $\boldsymbol{\Phi} = \mathbf{q}_k$ in (\ref{A_priori_pressure_4}) and using Young's inequality and (\ref{A_priori_pressure_5}), we obtain
\begin{equation*}
\normH{p_k}^2 = \mathcal{F}_k(\mathbf{q}_k)\leq \norm{\mathcal{F}_k}_{(\mathbf{H}^1)^*}\norm{\mathbf{q}_k}_{\mathbf{H}^1}\leq C_d\norm{\mathcal{F}_k}_{(\mathbf{H}^1)^*}\normH{p_k}\leq \frac{C_d^2}{2}\norm{\mathcal{F}_k}_{(\mathbf{H}^1)^*}^2 + \frac{1}{2}\normH{p_k}^2.
\end{equation*}
This implies
\begin{equation}
\label{A_priori_pressure_6}\normH{p_k}\leq C_d\norm{\mathcal{F}_k}_{(\mathbf{H}^1)^*}.
\end{equation}
Using Young's and Hölder's inequalities and (\ref{A_priori_pressure_3}), (\ref{A_priori_pressure_6}), we obtain

\begin{align}
\nonumber\intT \normH{p_k}^{\frac{4}{3}}\d t &\leq C\intT 1+\norm{\bv_k}_{\mathbf{H}^1}^{\frac{4}{3}} + \norm{\mu_k}_{L^3}^{\frac{4}{3}}\norm{\na\varphi_k}_{\mathbf{L}^2}^{\frac{4}{3}} + \norm{N_{\sigma,k}}_{L^3}^{\frac{4}{3}}\norm{\na\sigma_k}_{\mathbf{L}^2}^{\frac{4}{3}}\d t\\
\label{A_priori_pressure_7}&\leq C(1+\norm{\bv_k}_{L^2(\mathbf{H}^1)}^{\frac{4}{3}} + \norm{\mu_k}_{L^2(L^3)}^{\frac{4}{3}}\norm{\na\varphi_k}_{L^4(\mathbf{L}^2)}^{\frac{4}{3}} + \norm{N_{\sigma,k}}_{L^4(L^3)}^{\frac{4}{3}}\norm{\na\sigma_k}_{L^2(\mathbf{L}^2)}^{\frac{4}{3}})
\end{align}
Due to (\ref{central_estimate_1})-(\ref{central_estimate_2}) and (\ref{central_estimate_3})-(\ref{central_estimate_4}), the first three terms on the r.h.s. of (\ref{A_priori_pressure_7}) and the term $\norm{\na\sigma_k}_{L^2(\mathbf{L}^2)}^{\frac{4}{3}}$ are bounded. Thus, it remains to show that $N_{\sigma,k}\in L^4(L^3)$ with bounded norm. Using Gagliardo-Nirenberg with $j=0,~p=3,~m=1,~r=q=2$ yields
\begin{equation*}
\norm{N_{\sigma,k}}_{L^3}\leq C\normH{N_{\sigma,k}}^{\frac{1}{2}}\normV{N_{\sigma,k}}^{\frac{1}{2}}.
\end{equation*}
Therefore, it holds
\begin{equation*}
\intT \norm{N_{\sigma,k}}_{L^3}^4\d t \leq C\intT \normH{N_{\sigma,k}}^2\normV{N_{\sigma,k}}^2\leq C\norm{N_{\sigma,k}}_{L^{\infty}(L^2)}^2\norm{N_{\sigma,k}}_{L^2(H^1)}^2.
\end{equation*}
Due to the specific form of $N_{\sigma,k}$, applying (\ref{central_estimate_1})-(\ref{central_estimate_2}) and (\ref{central_estimate_3})-(\ref{central_estimate_4}) implies
\begin{equation*}
\norm{N_{\sigma,k}}_{L^4(L^3)}\leq C.
\end{equation*}
Consequently, from (\ref{A_priori_pressure_7}) we obtain
\begin{equation}
\label{A_priori_pressure_8}\norm{p_k}_{L^{\frac{4}{3}}(L^2)}\leq C,
\end{equation}
where $C$ is independent of $k\in\N$.

\subsection{Higher order estimates for $\varphi_k$}

In this section, we will show the inequality
\begin{equation}
\label{higher_order_estimates_1}\norm{\varphi_k}_{L^2(H^2)}\leq C.
\end{equation}
Using the Gagliardo-Nirenberg inequality and the continuous embedding $H^1\hookrightarrow L^6$, we have the following estimate:
\begin{equation}
\label{higher_order_estimates_2}\norm{\varphi_k}_{L^{\infty}}\leq C\norm{\varphi_k}_{H^1}^{\frac{1}{2}}\norm{\varphi_k}_{H^2}^{\frac{1}{2}}.
\end{equation}
Using elliptic regularity theory, this implies
\begin{equation}
\label{higher_order_estimates_3}\norm{\varphi_k}_{L^{\infty}}\leq C\norm{\varphi_k}_{H^1}^{\frac{1}{2}}(\normH{\varphi_k}^{\frac{1}{2}}+\norm{\Delta\varphi_k}^{\frac{1}{2}}).
\end{equation}
Chosing $v=\lambda_ja_j^kw_j$ in (\ref{approx_problem_eq_1b}), integrating by parts and summing the resulting equations over $j=1,...,k,$ yields
\begin{equation}
\label{higher_order_estimates_4}\epsilon\normH{\Delta\varphi_k}^2 = \intO \na\mu_k\cdot\na\varphi_k - \epsilon^{-1}\psi''(\varphi_k)|\na\varphi_k|^2 + \chi_{\varphi}\na\sigma_k\cdot\na\varphi_k\d x.
\end{equation}
Using Hölder's inequality and the assumption on $\psi$, we therefore get
\begin{equation}
\label{higher_order_estimates_5}\epsilon\normH{\Delta\varphi_k}^2\leq  \norm{\na\mu_k}_{\mathbf{L}^2}\norm{\na\varphi_k}_{\mathbf{L}^2} + \chi_{\varphi}\norm{\na\sigma_k}_{\mathbf{L}^2}\norm{\na\varphi_k}_{\mathbf{L}^2} + \intO C(1+|\varphi_k|^q)|\na\varphi_k|^2\d x.
\end{equation}
Integrating in time from $0$ to $T$, applying Hölder's inequality and (\ref{central_estimate_1}) gives
\begin{align}
\nonumber \epsilon\norm{\Delta\varphi_k}_{L^2(L^2)}^2 &\leq \norm{\na\mu_k}_{L^2(\mathbf{L}^2)}\norm{\na\varphi_k}_{L^2(\mathbf{L}^2)} + \chi_{\varphi}\norm{\na\sigma_k}_{L^2(\mathbf{L}^2)}\norm{\na\varphi_k}_{L^2(\mathbf{L}^2)}\\
\nonumber&\quad + C\intT\intO (1+|\varphi_k|^q)|\na\varphi_k|^2\d x\d t\\ 
\label{higher_order_estimates_6}&\leq C\left(1+ \intT\intO |\varphi_k|^q|\na\varphi_k|^2\d x\d t\right).
\end{align}
In the case $q=0$, applying (\ref{central_estimate_1}) gives
\begin{equation}
\label{higher_order_estimates_7}\intT\intO |\varphi_k|^q|\na\varphi_k|^2\d x\d t = \norm{\na\varphi_k}_{L^2(\mathbf{L}^2)}^2\leq C.
\end{equation}
In the case $q\in (0,4)$, we use Hölder's inequality and (\ref{higher_order_estimates_3}) to calculate
\begin{align}
\nonumber \intT\intO |\varphi_k|^q|\na\varphi_k|^2\d x\d t&\leq C\intT \norm{\varphi_k}_{L^{\infty}}^q\norm{\na\varphi_k}_{\mathbf{L}^2}^2\d t\\
\nonumber &\leq C\intT \normV{\varphi_k}^2\normV{\varphi_k}^{\frac{q}{2}}(\normH{\varphi_k}^{\frac{q}{2}} + \normH{\Delta\varphi_k}^{\frac{q}{2}})\d t\\
\label{higher_order_estimates_8}&\leq C \left(\norm{\varphi_k}_{L^{\infty}(H^1)}^{q+2} + \intT\normV{\varphi_k}^{\frac{q+4}{2}}\normH{\Delta\varphi_k}^{\frac{q}{2}}\d t\right).
\end{align}
Observing $\frac{4}{q}>1$, we can use Young's generalised inequality to estimate the last integral on the r.h.s. of (\ref{higher_order_estimates_8}) by
\begin{equation}
\label{higher_order_estimates_9}C\intT\normV{\varphi_k}^{\frac{q+4}{2}}\normH{\Delta\varphi_k}^{\frac{q}{2}}\d t \leq C\norm{\varphi_k}_{L^{\infty}(H^1)}^{\frac{2(q+4)}{4-q}} + \frac{\epsilon}{2}\norm{\Delta\varphi_k}_{L^2(L^2)}^2.
\end{equation}
Consequently, from (\ref{higher_order_estimates_6})-(\ref{higher_order_estimates_9}) we obtain
\begin{equation}
\frac{\epsilon}{2}\norm{\Delta\varphi_k}_{L^2(L^2)}^2\leq C.
\end{equation}
Using elliptic regularity theory and (\ref{central_estimate_2}), this implies 
\begin{equation}
\label{central_estimate_5}\norm{\varphi_k}_{L^2(H^2)}\leq C.
\end{equation}
Summarizing the estimates (\ref{central_estimate_1})-(\ref{central_estimate_2}), (\ref{central_estimate_3})-(\ref{central_estimate_4}) and (\ref{central_estimate_5}), we deduce that
\begin{equation}
\label{central_estimate_6}\norm{\varphi_k}_{L^{\infty}(H^1)\cap L^2(H^2)} + \norm{\sigma_k}_{L^{\infty}(L^2)\cap L^2(H^1)} + \norm{\mu_k}_{L^2(H^1)} + \norm{\bv_k}_{L^2(\mathbf{H}^1)}\leq C.
\end{equation}

\subsection{Regularity for the convection terms and the time derivatives}
By Hölder's inequality and the continuous embedding $\mathbf{H}^1\hookrightarrow\mathbf{L}^6$, we observe that
\begin{align*}
\norm{\na\varphi_k\cdot\bv_k}_{L^2(0,T;L^{\frac{3}{2}})}^2 &= \intT \norm{\na\varphi_k\cdot\bv_k}_{L^{\frac{3}{2}}}^2\d t\leq \intT \norm{\bv_k}_{\mathbf{L}^6}^2\norm{\na\varphi_k}_{\mathbf{L}^2}^2\d t\\
&\quad\leq C\intT \norm{\bv_k}_{\mathbf{H}^1}^2\normV{\varphi_k}^2\d t\\
&\quad\leq C\norm{\varphi_k}_{L^{\infty}(0,T;H^1)}^2\norm{\bv_k}_{L^2(0,T;\mathbf{H}^1)}^2.
\end{align*}
Using the boundedness of $\Gamma_{\bv}$ and (\ref{central_estimate_2}), (\ref{central_estimate_4}), we see that
\begin{equation*}
\norm{\varphi_k\Gamma_{\bv,k}}_{L^2(0,T;L^{\frac{3}{2}})}^2 \leq \gamma_0^2\norm{\varphi_k}_{L^2(0,T;L^{\frac{3}{2}})}^2\leq C,
\end{equation*}
with a constant $C$ depending on $\gamma_0$.
From the last two inequalities and (\ref{central_estimate_1})-(\ref{central_estimate_2}), (\ref{central_estimate_3})-(\ref{central_estimate_4}), we deduce that
\begin{equation}
\label{Apriori_convection_terms_1}\norm{\div(\varphi_k\bv_k)}_{L^2(0,T;L^{\frac{3}{2}})}\leq C.
\end{equation}
Taking an arbitrary $\zeta\in L^4(0,T;H^1)$ with coefficients $\{\zeta_{kj}\}_{1\leq j\leq k}$ such that $\mathbb{P}_k\zeta = \sum_{j=1}^{k}\zeta_{kj}w_j$ and integrating by parts, we obtain
\begin{equation}
\label{proof_convection_terms_1}\intT\intO \div(\sigma_k\bv_k)\mathbb{P}_k\zeta\d x\d t = \intT\int_{\pa\Omega}\mathbb{P}_k\zeta\sigma_k\bv_k\cdot\mathbf{n}\d \mathcal{H}^{d-1}\d t - \intT\intO \sigma_k\bv_k\cdot\na\mathbb{P}_k\zeta\d x\d t.
\end{equation}
By the Gagliardo-Nirenberg inequality (\ref{lemma_Gagliardo_Nirenberg}) with $j=0,~p=3,~m=1,~r=2,~q=2$, we have
\begin{equation*}
\norm{\sigma_k}_{L^3}\leq C \normH{\sigma_k}^{\frac{1}{2}}\normV{\sigma_k}^{\frac{1}{2}}.
\end{equation*}
Then, by Hölder's inequality and the continuous embedding $\mathbf{H}^1\hookrightarrow\mathbf{L}^6$, we can estimate the second term on the r.h.s. of (\ref{proof_convection_terms_1}) by
\begin{align}
\nonumber\left|\intT\intO \sigma_k\bv_k\cdot \na \mathbb{P}_k\zeta\d x\d t\right|&\leq \intT \norm{\sigma_k}_{L^3}\norm{\bv_k}_{\mathbf{L}^6}\norm{\na \mathbb{P}_k\zeta}_{\mathbf{L}^2}\d t\\
\nonumber&\leq C\intT \normH{\sigma_k}^{\frac{1}{2}}\normV{\sigma_k}^{\frac{1}{2}}\norm{\bv_k}_{\mathbf{H}^1}\normV{\zeta}\d t\\
\label{proof_convection_terms_2}&\leq C \norm{\sigma_k}_{L^{\infty}(L^2)}^{\frac{1}{2}}\norm{\sigma_k}_{L^2(H^1)}^{\frac{1}{2}}\norm{\bv_k}_{L^2(\mathbf{H}^1)}\norm{\zeta}_{L^4(H^1)}
\end{align}
Furthermore, using (\ref{interpolation_inequality_boundary}) with $r=q=2$ (hence $\alpha = 0$) gives
\begin{equation*}
\norm{\sigma_k}_{L^2(\pa\Omega)}\leq C(\norm{\sigma_k}_{L^2} + \norm{\sigma_k}_{L^2}^{\frac{1}{2}}\norm{\sigma_k}_{H^1}^{\frac{1}{2}})
\end{equation*}
Using (\ref{central_estimate_1}) and (\ref{central_estimate_3}), this implies
\begin{equation}
\label{proof_convection_terms_3}\norm{\sigma_k}_{L^4(L^2(\pa\Omega))}\leq C.
\end{equation}
Now, using Hölder's inequality and the trace theorem, we obtain
\begin{align}
\nonumber \left|\intT\int_{\pa\Omega}\mathbb{P}_k\zeta\sigma_k\bv_k\cdot\mathbf{n}\d \mathcal{H}^{d-1}\d t\right|&\leq \intT \norm{\sigma_k}_{L^2(\pa\Omega)}\norm{\bv_k}_{\mathbf{L}^4(\pa\Omega)}\norm{\mathbb{P}_k\zeta}_{L^4(\pa\Omega)}\d t\\
\nonumber &\leq C\intT \norm{\sigma_k}_{L^2(\pa\Omega)}\norm{\bv_k}_{\mathbf{H}^1}\norm{\zeta}_{H^1}\d t\\
\label{proof_convection_terms_4}&\leq C\norm{\sigma_k}_{L^4(L^2(\pa\Omega))}\norm{\bv_k}_{L^2(\mathbf{H}^1)}\norm{\zeta}_{L^4(H^1)}.
\end{align}
Hence, from (\ref{central_estimate_1}), (\ref{central_estimate_3}) and (\ref{proof_convection_terms_1})-(\ref{proof_convection_terms_4}) we get
\begin{equation}
\label{Apriori_convections_terms_2}\norm{\div(\sigma_k\bv_k)}_{L^{\frac{4}{3}}(0,T;(H^1)^*)}\leq C.
\end{equation}
Now, taking $v=\zeta_{kj}w_j$ in (\ref{approx_problem_eq_1c}) and summing over $j=1,...,k$, integrating in time from $0$ to $T$ yields
\begin{align*}
\left|\intT\intO \pa_t\sigma_k\zeta\d x\d t\right| &\leq \intT n_1(\chi_{\sigma}\norm{\na\sigma_k}_{\mathbf{L}^2}+\chi_{\varphi}\norm{\na\varphi_k}_{\mathbf{L}^2})\norm{\na \mathbb{P}_k\zeta}_{\mathbf{L}^2} + \normH{\Gamma_{\sigma,k}}\norm{\mathbb{P}_k\zeta}_{L^2}\d t\\
&\quad + \intT \norm{\div(\sigma_k\bv_k)}_{(H^1)^*}\norm{\mathbb{P}_k\zeta}_{H^1}\d t \\
&\quad+ \intT b(\norm{\sigma_{\infty}}_{L^2(\pa\Omega)}+\norm{\sigma_k}_{L^2(\pa\Omega)})\norm{\mathbb{P}_k\zeta}_{L^2(\pa\Omega)}\d t.
\end{align*}
Using (\ref{central_estimate_1})-(\ref{central_estimate_2}), (\ref{central_estimate_3})-(\ref{central_estimate_4}), we obtain
\begin{equation*}
\norm{\Gamma_{\sigma,k}}_{L^2(L^2)} \leq C(R_0,|\Omega|,T)(1 + \norm{\varphi_k}_{L^2(L^2)}+\norm{\sigma_k}_{L^2(L^2)}+\norm{\mu_k}_{L^2(L^2)})\leq C.
\end{equation*}
Then, Hölder's inequality and the trace theorem yields
\begin{equation*}
\left|\intT\intO \pa_t\sigma_k\zeta\d x\d t\right|\leq C(1+\norm{\div(\sigma_k\bv_k)}_{L^{\frac{4}{3}}((H^1)^*)})\norm{\zeta}_{L^4(H^1)}.
\end{equation*}
By taking the supremum over all $\zeta\in L^4(H^1)$ and using (\ref{central_estimate_1})-(\ref{central_estimate_2}), (\ref{central_estimate_3})-(\ref{central_estimate_4}) and (\ref{Apriori_convections_terms_2}), we end up with
\begin{equation}
\label{Apriori_time_derivatives_1}\norm{\pa_t\sigma_k}_{L^{\frac{4}{3}}((H^1)^*)}\leq C.
\end{equation}
With similar arguments, we can show that
\begin{equation}
\label{Apriori_time_derivatives_2}\norm{\pa_t\varphi_k}_{L^2((H^1)^*)}.
\end{equation}
Notice that we have lower time regularity for the time derivative of $\varphi_k$ compared to the convection term since the regularity of the time derivative depends on the term $\na\mu_k$.

\section{Passing to the limit}
At this point, we summarise the estimates (\ref{A_priori_pressure_8}), (\ref{central_estimate_6})-(\ref{Apriori_convection_terms_1}), (\ref{Apriori_convections_terms_2})-(\ref{Apriori_time_derivatives_2}) to deduce

\begin{align}
\nonumber&\norm{\varphi_k}_{L^{\infty}(H^1)\cap L^2(H^2)\cap W^{1,2}((H^1)^*)} + \norm{\sigma_k}_{L^{\infty}(L^2)\cap L^2(H^1)\cap W^{1,\frac{4}{3}}((H^1)^*)} + \norm{\mu_k}_{L^2(H^1)}\\
\label{central_estimate}&\quad + \norm{\div(\varphi_k\bv_k)}_{L^2(L^{\frac{3}{2}})} + \norm{\div(\sigma_k\bv_k)}_{L^{\frac{4}{3}}((H^1)^*)} + \norm{\bv_k}_{L^2(\mathbf{H}^1)}+ \norm{p_k}_{L^{\frac{4}{3}}(L^2)}\leq C.
\end{align}
Using  standard compactness arguments (Aubin-Lions theorem (see \cite[Sec. 8, Cor. 4]{Simon}) and reflexive weak compactness), the compact embeddings
\begin{equation*}
H^{j+1} = W^{j+1,2}\hookrightarrow\hookrightarrow W^{j,r}\quad \forall j\in \N_0,~1\leq r<6,
\end{equation*}
and $L^2\hookrightarrow\hookrightarrow (H^1)^*$, $H^1\hookrightarrow\hookrightarrow H^{\frac{1}{2}}$, we obtain, at least for a subsequence which will again be labelled by $k$, the following convergence results:
\begin{alignat*}{3}
\varphi_k&\ra\varphi&&\quad\text{weakly-}*&&\quad\text{ in }L^{\infty}(H^1)\cap L^2(H^2)\cap W^{1,2}((H^1)^*),\\
\sigma_k&\ra \sigma &&\quad\text{weakly-}*&&\quad\text{ in }L^{\infty}(L^2)\cap L^2(H^1)\cap W^{1,\frac{4}{3}}((H^1)^*),\\
\mu_k&\ra \mu&&\quad\text{weakly}&&\quad\text{ in }L^2(H^1),\\
p_k&\ra p&&\quad\text{weakly}&&\quad\text{ in }L^{\frac{4}{3}}(L^2),\\
\bv_k&\ra\bv&&\quad\text{weakly}&&\quad\text{ in }L^2(\mathbf{H}^1),\\
\div(\varphi_k\bv_k)&\ra \tau&&\quad\text{weakly}&&\quad\text{ in }L^2(L^{\frac{3}{2}}),\\
\div(\sigma_k\bv_k)&\ra \theta&&\quad\text{weakly}&&\quad\text{ in }L^{\frac{4}{3}}((H^1)^*),\\
\div(\bv_k)&\ra \div(\bv)&&\quad\text{weakly}&&\quad \text{ in }L^2(L^2),
\end{alignat*}
for some limit functions $\tau\in L^2(L^{\frac{3}{2}})$ and $\theta\in L^{\frac{4}{3}}((H^1)^*)$. Furthermore, we have the strong convergences
\begin{alignat*}{3}
\varphi_k&\ra\varphi&&\quad\text{strongly}&&\quad\text{ in }C^0(L^r)\cap L^2(W^{1,r})\text{ and a.e. in }Q,\\
\sigma_k&\ra\sigma&&\quad\text{strongly}&&\quad\text{ in }C^0((H^1)^*)\cap L^2(L^r)\cap L^2(H^{\frac{1}{2}})\text{ and a.e. in }Q,
\end{alignat*}
for $r\in [1,6)$. From now on, we fix $1\leq j\leq k$ and $\xi\in L^2,~\boldsymbol{\Phi}\in \mathbf{H}^1$, $\delta\in C_0^{\infty}(0,T)$. Then, since the eigenfunctions $\{w_j\}_{j\in\N}$ belong to $H^2$, we observe that $\delta w_j\in C^{\infty}(H^2)$ for all $j\in \N$. Furthermore, we have $\delta\xi\in C^{\infty}(L^2),~\delta\boldsymbol{\Phi}\in C^{\infty}(\mathbf{H}^1)$. Inserting $v=w_j$ in (\ref{approx_problem_eq_1a})-(\ref{approx_problem_eq_1c}), multiplying the resulting equations with $\delta$ and integrating over $(0,T)$ yields
\begin{align}
\label{limiting_equations_1}&\intT \delta(t)\left( \intO (\pa_t\varphi_k - \Gamma_{\varphi,k} + \na\varphi_k\cdot\bv_k + \varphi_k\Gamma_{\bv,k})w_j + m(\varphi_k)\na\mu_k\cdot\na w_j\d x\right)\d t = 0,\\
\label{limiting_equations_2}&\intT \delta(t)\left(\intO (\mu_k - \epsilon^{-1}\psi'(\varphi_k)+ \chi_{\varphi}\sigma_k)w_j - \epsilon\na\varphi_k\cdot\na w_j\d x\right)\d t =0,\\
\nonumber&\intT\delta(t)\left(\intO (\pa_t\sigma_k + \Gamma_{\sigma,k} + \na\sigma_k\cdot\bv_k + \sigma_k\Gamma_{\bv,k})w_j +n(\varphi_k)\na N_{\sigma,k}\na w_j\d x\right)\d t\\
\label{limiting_equations_3}&\quad - \intT\delta(t)\left(\int_{\pa\Omega}b(\sigma_{\infty}-\sigma_k)w_j\d\mathcal{H}^{d-1}\right)\d t =0.
\end{align}
Furthermore, we take the $\mathbf{L}^2$-scalar product of (\ref{approx_problem_eq_1d}) with $\boldsymbol{\Phi}$, multiply with $\delta$ and integrate from $0$ to $T$ to obtain
\begin{equation}
\label{limiting_equations_4}\intT \delta(t)\left(\intO T(\bv_k,p_k)\colon \na\boldsymbol{\Phi} + \nu\bv_k\cdot \boldsymbol{\Phi} - (\mu_k\na\varphi_k + N_{\sigma,k}\na\sigma_k)\cdot\boldsymbol{\Phi}\d x\right)\d t =0,
\end{equation}
where we used (\ref{approx_problem_eq_1f}). With similar arguments, (\ref{approx_problem_eq_1e}) gives
\begin{equation}
\label{limiting_equations_5}\intT\delta(t)\left( \intO \div(\bv_k)\xi\d x\right)\d t = \intT \delta(t)\left(\intO \Gamma_{\bv,k}\xi\d x\right)\d t.
\end{equation}
Now, we want to pass to the limit in (\ref{limiting_equations_1})-(\ref{limiting_equations_5}).\newline\newline
\underline{\textbf{Step 1:} (\ref{limiting_equations_1})}
Since $\delta w_j\in C^{\infty}(H^2)\hookrightarrow L^2((H^1)^*)$, we obtain
\begin{equation}
\label{limit_1}\intT \intO \pa_t\varphi_k\delta w_j\d x\d t \ra \intT \delta(t)\langle \pa_t\varphi{,}w_j\rangle_{H^1}\d t\quad\text{ as }k\ra\infty.
\end{equation}
By continuity of $m(\cdot)$ and since $\varphi_k\ra\varphi$ a.e. in $Q$ as $k\ra\infty$, we observe that $m(\varphi_k)\ra m(\varphi)$ a.e. in $Q$. Using the boundedness of $m(\cdot)$ and applying Lebesgue dominated convergence theorem to $(m(\varphi_k)-m(\varphi))^2|\delta|^2|\na w_j|^2$, we obtain
\begin{equation*}
\norm{(m(\varphi_k)-m(\varphi))\delta\na w_j}_{L^2(Q)}\ra 0\quad\text{ as }k\ra\infty.
\end{equation*}
Then, by weak convergence $\na\mu_k\rightharpoonup\na\mu$ in $L^2(Q)$ as $k\ra\infty$, we have
\begin{equation}
\label{limit_2}\intT\intO \delta m(\varphi_k)\na w_j\cdot\na\mu_k\d x\d t\ra \intT\intO \delta m(\varphi)\na w_j\cdot\na\mu\d x\d t\quad\text{ as }k\ra\infty.
\end{equation}
Using the continuous embedding $H^2\hookrightarrow L^{\infty}$, we have
\begin{align*}
\intT \intO |\delta|^2 |w_j|^2 |\na\varphi_k-\na\varphi|^2\d x\d t &\leq \intT |\delta|^2\norm{\na\varphi_k-\na\varphi}_{\mathbf{L}^2}^2\norm{w_j}_{L^{\infty}}^2\d t\\
&\leq C\norm{\delta}_{L^{\infty}(0,T)}^2\norm{w_j}_{H^2}^2\norm{\varphi_k-\varphi}_{L^2(H^1)}^2\\
&\ra 0\quad \text{ as }k\ra\infty.
\end{align*}
Therefore, $\delta w_j\na\varphi_k\ra\delta w_j\na\varphi$ strongly in $L^2(\mathbf{L}^2)$ as $k\ra\infty$. Then, by the product of weak-strong convergence, we obtain
\begin{equation}
\label{limit_3}\intT \intO \delta w_j\na\varphi_k\cdot\bv_k\d x\d t\ra \intT\intO \delta w_j\na\varphi\cdot\bv \d x\d t\quad\text{ as }k\ra\infty.
\end{equation}
Using the boundedness and continuity of $\Gamma_{\bv}(\cdot{,}\cdot)$, with similar arguments as for (\ref{limit_2}) we obtain
\begin{equation*}
\norm{(\Gamma_{\bv}(\varphi_k,\sigma_k)-\Gamma_{\bv}(\varphi,\sigma))\delta w_j}_{L^2(Q)}\ra 0\quad\text{ as }k\ra\infty,
\end{equation*}
which implies
\begin{equation}
\label{limit_7}\intT\intO \delta w_j\varphi_k\Gamma_{\bv}(\varphi_k,\sigma_k)\d x\d t\ra \intT\intO \delta w_j\varphi\Gamma_{\bv}(\varphi,\sigma)\d x\d t\quad\text{ as }k\ra\infty.
\end{equation}
In particular, from (\ref{limit_3}) and (\ref{limit_7}) we can conclude that $\div(\varphi\bv) = \tau$.
Now, we recall the specific form of $N_{\varphi,k}$ given by $N_{\varphi,k} = \Lambda_{\varphi}(\varphi_k,\sigma_k) - \theta_{\varphi}(\varphi_k,\sigma_k)\mu_k$. Using that $\varphi_k\ra\varphi$ and $\sigma_k\ra\sigma$ a.e. in $Q$ together with the continuity and boundedness of $\theta_{\varphi}(\cdot{,}\cdot)$, Lebesgue's dominated convergence theorem implies
\begin{equation*}
\intT \intO |\delta w_j(\theta_{\varphi}(\varphi_k,\sigma_k)-\theta_{\varphi}(\varphi,\sigma))|^2\d x\d t\ra 0\quad\text{ as }k\ra\infty.
\end{equation*}
Therefore, $\delta w_j\theta_{\varphi}(\varphi_k,\sigma_k)\ra \delta w_j\theta_{\varphi}(\varphi,\sigma)$ strongly in $L^2(Q)$ as $k\ra\infty$. Together with the weak convergence $\mu_k\rightharpoonup\mu$ in $L^2(Q)$, we conclude that
\begin{equation}
\label{limit_4}\intT\intO \delta w_j\theta_{\varphi}(\varphi_k,\sigma_k)\mu_k\d x\d t\ra \intT\intO \delta w_j\theta_{\varphi}(\varphi,\sigma)\mu\d x\d t\quad\text{ as }k\ra\infty.
\end{equation}
We now analyse the other term in the definition of $\Gamma_{\varphi,k}$. 
Applying the inequality $||a|-|b||\leq |a-b|$, we obtain
\begin{equation*}
\intT \intO |(\delta w_j)(|\varphi_k|-|\varphi|)|\d x\d t\leq \norm{\delta w_j}_{L^2(Q)}\norm{\varphi_k-\varphi}_{L^2(Q)}\xrightarrow{k\ra\infty}0
\end{equation*}
and 
\begin{equation*}
\intT \intO |(\delta w_j)(|\sigma_k|-|\sigma|)|\d x\d t\leq \norm{\delta w_j}_{L^2(Q)}\norm{\sigma_k-\sigma}_{L^2(Q)}\xrightarrow{k\ra\infty}0.
\end{equation*}
This implies
\begin{equation*}
R_0(1+|\varphi_k|+|\sigma_k|)|\delta w_j|\ra R_0(1+|\varphi|+|\sigma|)|\delta w_j|\quad\text{ strongly in }L^1(Q)\text{ as }k\ra\infty.
\end{equation*}
Since $\varphi_k\ra\varphi$ and $\sigma_k\ra\sigma$ a.e. in $Q$ as $k\ra\infty$, the continuity of $\Lambda(\cdot{,}\cdot)$ yields
\begin{equation*}
\delta w_j\Lambda(\varphi_k,\sigma_k)\ra \delta w_j\Lambda(\varphi,\sigma) \quad\text{ a.e. in }Q \text{ as }k\ra\infty.
\end{equation*}
Using
\begin{equation*}
|\delta w_j\Lambda(\varphi_k,\sigma_k)|\leq |\delta w_j|R_0(1+|\varphi_k|+|\varphi_k|)\in L^1(\Omega\times (0,T))\quad \forall k\geq 1,
\end{equation*}
by the generalised Lebesgue dominated convergence theorem (see \cite[3.25, p.60]{Alt}) we obtain
\begin{equation*}
\intT\intO \delta w_j \Lambda_{\varphi}(\varphi_k,\sigma_k)\d x\d t\ra \intT \intO \delta w_j\Lambda_{\varphi}(\varphi,\sigma)\d x\d t\quad\text{ as }k\ra\infty.
\end{equation*}
Together with (\ref{limit_4}), this implies
\begin{equation}
\label{limit_5}\intT\intO \delta w_j \Gamma_{\varphi}(\varphi_k,\sigma_k,\mu_k)\d x\d t\ra \intT \intO \delta w_j\Gamma_{\varphi}(\varphi,\sigma,\mu)\d x\d t\quad\text{ as }k\ra\infty.
\end{equation}
\underline{\textbf{Step 2}:} We now want to analyse (\ref{limiting_equations_2}). Since $\mu_k\rightharpoonup\mu,~\sigma_k\rightharpoonup\sigma$ and $\na\varphi_k\rightharpoonup\na\varphi$ in $L^2(Q)$ and $L^2(\mathbf{L}^2)$, we easily deduce
\begin{align}
\nonumber &\intT\delta(t)\left(\intO (\mu_k + \chi_{\varphi}\sigma_k)w_j -\epsilon\na\varphi_k\cdot\na w_j\d x\right)\d t\\
\label{limit_8}&\ra \intT\delta(t)\left(\intO (\mu + \chi_{\varphi}\sigma)w_j -\epsilon\na\varphi\cdot\na w_j\d x\right)\d t\quad \text{ as }k\ra\infty.
\end{align}
If the derivative $\psi'(\cdot)$ satisfies the linear growth condition (\ref{assumptions_on_psi_2}), we can use similar arguments as for (\ref{limit_5}) to deduce that 
\begin{equation}
\intT\intO \epsilon^{-1}\psi'(\varphi_k)\delta w_j\d x\d t\ra \intT\intO \epsilon^{-1}\psi'(\varphi)\delta w_j\d x\d t\quad\text{ as }k\ra\infty.
\end{equation}
For potentials satisfying (\ref{assumption_on_psi_4}), we refer to the argument in \cite[§3.1.2]{GarckeLam2}.\newline
\underline{\textbf{Step 3}:} We now want to pass to the limit in (\ref{limiting_equations_5}). Since $\varphi_k\ra\varphi,~\sigma_k\ra\sigma$ a.~e. in $Q$ as $k\ra\infty$, the continuity and boundedness of $\Gamma_{\bv}$ and similar arguments as for (\ref{limit_2}) imply
\begin{equation*}
\intT\intO \delta(t)\Gamma_{\bv,k}\xi\d x\d t\ra \intT\intO \delta(t)\Gamma_{\bv}(\varphi,\sigma)\xi\d x\d t\quad \text{ as }k\ra\infty.
\end{equation*}
Recalling the weak convergence $\div(\bv_k)\ra \div(\bv)$ in $L^2(L^2)$ as $k\ra\infty$, we deduce
\begin{equation*}
\intT\intO \delta(t)\div(\bv_k)\xi\d x\d t\ra \intT\intO \delta(t)\div(\bv)\xi\d x\d t\quad \text{ as }k\ra\infty.
\end{equation*}
This allows us to pass to the limit $k\ra\infty$ in (\ref{limiting_equations_5}) to obtain
\begin{equation}
\label{limit_14}\intT\delta(t)\intO \div(\bv)\xi\d x\d t = \intT\delta(t)\intO \Gamma_{\bv}(\varphi,\sigma)\xi\d x\d t.
\end{equation}
In particular, since this holds for all $\delta\in C_0^{\infty}(0,T)$ and all $\xi\in L^2$, we have
\begin{equation}
\label{limit_14a}\div(\bv) = \Gamma_{\bv}(\varphi,\sigma)\quad\text{a.~e. in }Q.
\end{equation}
\underline{\textbf{Step 4}:}
With similar arguments as for (\ref{limit_1})-(\ref{limit_2}) and (\ref{limit_5}), we obtain
\begin{align}
\label{limit_11}\intT\intO \pa_t\sigma_k\delta w_j\d x\d t &\ra \intT\delta(t)\langle \pa_t\sigma{,}w_j\rangle_{H^1,(H^1)^*}\d t,\\
\label{limit_12}\intT\intO \delta n(\varphi_k)\na N_{\sigma,k}\cdot\na w_j\d x\d t &\ra \intT\intO \delta n(\varphi)(\chi_{\sigma}\na\sigma-\chi_{\varphi}\na\varphi)\cdot\na w_j\d x\d t,\\
\label{limit_6}\intT\intO \delta w_j\Gamma_{\sigma}(\varphi_k,\sigma_k,\mu_k)\d x\d t&\ra \intT\intO \delta w_j\Gamma_{\sigma}(\varphi,\sigma,\mu)\d x\d t
\end{align}
as $k\ra\infty$. For the boundary term in (\ref{limiting_equations_3}), we first recall the continuous embedding $H^1(\Omega)\hookrightarrow L^4(\pa\Omega)$. Then, by the weak convergence of $\sigma_k\rightharpoonup\sigma$ in $L^2(\Sigma)$, we conclude
\begin{equation}
\label{limit_13}\intT \delta(t)\left(\int_{\pa\Omega}\sigma_k w_j\d \mathcal{H}^{d-1}\right)\d t\ra \intT \delta(t)\left(\int_{\pa\Omega}\sigma w_j\d \mathcal{H}^{d-1}\right)\d t\quad\text{ as }k\ra\infty.
\end{equation}
To pass to the limit in the convection term of (\ref{approx_problem_eq_1c}), we first want to show that
\begin{equation}
\label{limit_16}\intT\intO \sigma\Gamma_{\bv}(\varphi,\sigma)\delta w_j \d x\d t = \intT\intO \sigma\div(\bv)\delta w_j\d x\d t.
\end{equation}
Indeed, a short calculation yields
\begin{align*}
\intT\intO |\delta|^2 |w_j|^2|\sigma_k-\sigma|^2\d x\d t &\leq \intT |\delta|^2\norm{w_j}_{L^6}^2\norm{\sigma_k-\sigma}_{L^3}^2\d t\\
&\leq C\norm{\delta}_{L^{\infty}(0,T)}^2\normV{w_j}^2\norm{\sigma_k-\sigma}_{L^2(L^3)}^2\\
&\ra 0\quad\text{ as }k\ra\infty,
\end{align*}
where we used that $\sigma_k\ra\sigma$ strongly in $L^2(L^3)$. Therefore, we obtain that $\sigma_k\delta w_j\ra \sigma\delta w_j$ strongly in $L^2(L^2)$. With similar arguments as for (\ref{limit_14}), this implies (\ref{limit_16}).
Now, as $\delta w_j\in C^{\infty}(H^2)\hookrightarrow L^4(H^1)$, the weak convergence $\div(\sigma_k\bv_k)\rightharpoonup \theta$ in $L^{\frac{4}{3}}((H^1)^*)$ implies
\begin{equation}
\label{limit_17}\intT\intO \div(\sigma_k\bv_k)\delta w_j\d x\d t\ra \intT \delta(t)\langle \theta{,}w_j\rangle_{H^1,(H^1)^*}\d t\quad\text{ as }k\ra\infty.
\end{equation}
Using integration by parts, we see that
\begin{equation}
\label{limit_16a}\intT\intO \div(\sigma_k\bv_k)\delta w_j\d x\d t = \intT\int_{\pa\Omega}\delta w_j\sigma_k\bv_k\cdot\mathbf{n}\d \mathcal{H}^{d-1}\d t- \intT\intO \delta\sigma_k\bv_k\cdot\na w_j\d x\d t.
\end{equation}
Recalling that $\sigma_k\ra\sigma$ strongly in $L^2(H^{\frac{1}{2}})$ and using the continuous embeddings $H^{\frac{1}{2}}\hookrightarrow L^2(\pa\Omega)$, $\mathbf{H}^1\hookrightarrow \mathbf{L}^4(\pa\Omega)$ resulting from the trace theorem, we have
\begin{alignat*}{3}
\sigma_k&\ra \sigma&&\quad\text{ strongly}&&\quad\text{ in }L^2(L^2(\pa\Omega)),\\
\bv_k&\ra\bv&&\quad\text{ weakly }&&\quad\text{ in }L^2(\mathbf{L}^4(\pa\Omega)).
\end{alignat*}
Again by the trace theorem and the continuous embeddings $H^2\hookrightarrow W^{1,6},~W^{\frac{5}{6},6}(\pa\Omega)\hookrightarrow L^{\infty}(\pa \Omega)$, we observe that $w_j\in H^2(\Omega)\hookrightarrow L^{\infty}(\pa\Omega)$. Since the outer unit normal $\mathbf{n}$ is continuous, we calculate
\begin{align*}
\intT\int_{\pa\Omega}|\delta|^2 |\mathbf{n}|^2|\sigma_k-\sigma|^2|w_j|^2\d  \mathcal{H}^{d-1}\d t &\leq \intT |\delta|^2 \norm{w_j}_{L^{\infty}(\pa\Omega)}^2\norm{\sigma_k-\sigma}_{L^2(\pa\Omega)}^2\d t\\
&\leq C\norm{\delta}_{L^{\infty}(0,T)}^2\norm{w_j}_{H^2}^2\norm{\sigma_k-\sigma}_{L^2(L^2(\pa\Omega))}^2\\
&\ra 0\quad\text{ as }k\ra\infty,
\end{align*}
meaning $\delta w_j\sigma_k\mathbf{n}\ra \delta w_j\sigma\mathbf{n}$ strongly in $L^2(\mathbf{L}^2(\pa\Omega))$ as $k\ra\infty$. Then, by the product of weak-strong convergence, we obtain
\begin{equation}
\label{limit_18}\intT\int_{\pa\Omega}\delta w_j\sigma_k\bv_k\cdot\mathbf{n}\d \mathcal{H}^{d-1}\d t \ra \intT\int_{\pa\Omega}\delta w_j\sigma\bv\cdot\mathbf{n}\d \mathcal{H}^{d-1}\d t\quad\text{ as }k\ra\infty.
\end{equation}
Furthermore, since $\sigma_k\ra\sigma$ strongly in $L^2(L^3)$, we get
\begin{align*}
\intT\intO |\delta|^2|\na w_j|^2|\sigma_k-\sigma|^2\d x\d t&\leq \intT |\delta|^2\norm{\na w_j}_{L^6}^2\norm{\sigma_k-\sigma}_{L^3}^2\d t\\
&\leq C\norm{\delta}_{L^{\infty}(0,T)}^2\norm{w_j}_{H^2}^2\norm{\sigma_k-\sigma}_{L^2(L^3)}\\
&\ra 0\quad\text{ as }k\ra\infty.
\end{align*}
Then, since $\bv_k\rightharpoonup \bv$ in $L^2(H^1)$, by the product of weak-strong convergence we have
\begin{equation}
\label{limit_19}\intT\intO \delta\sigma_k\bv_k\cdot\na w_j\d x\d t\ra \intT\intO \delta\sigma\bv\cdot\na w_j\d x\d t\quad\text{ as }k\ra\infty.
\end{equation}
Passing to the limit in (\ref{limit_17}) and using (\ref{limit_16a}), (\ref{limit_18})-(\ref{limit_19}), we obtain
\begin{equation}
\label{limit_20}\intT \delta(t)\langle\theta{,}w_j\rangle_{H^1,(H^1)^*}\d t = \intT\int_{\pa\Omega}\delta w_j\sigma\bv\cdot\mathbf{n}\d \mathcal{H}^{d-1}\d t - \intT\intO \delta\sigma\bv\cdot\na w_j\d x\d t.
\end{equation}
Again integrating by parts yields
\begin{equation}
\label{limit_21}\intT \delta(t)\langle\theta{,}w_j\rangle_{H^1,(H^1)^*}\d t = \intT\intO \div(\sigma\bv)\delta w_j\d x\d t,
\end{equation}
hence $\div(\sigma\bv) = \theta$ in the sense of distributions. In particular, by (\ref{limit_16}) we have
\begin{equation}
\label{limit_22}\intT \delta(t)\langle\theta{,}w_j\rangle_{H^1}\d t = \intT\intO \na\sigma\cdot\bv \delta w_j + \sigma\Gamma_{\bv}(\varphi,\sigma)\delta w_j\d x\d t.
\end{equation}
\underline{\textbf{Step 5}:}
Finally, we want to pass to the limit in (\ref{limiting_equations_4}). First of all, we recall that $\delta\boldsymbol{\Phi}\in C^{\infty}(\mathbf{H}^1)$. Then, by continuity of $\eta(\cdot)$, $\lambda(\cdot)$ and since $\varphi_k\ra\varphi$ a.e. in $Q$ as $k\ra\infty$, we observe that $\eta(\varphi_k)\ra \eta(\varphi)$, $\lambda(\varphi_k)\ra\lambda(\varphi)$ a.e. in $Q$. Using the boundedness of $\eta(\cdot)$ and $\lambda(\cdot)$, applying Lebesgue dominated convergence theorem to $(\eta(\varphi_k)-\eta(\varphi))^2|\delta|^2|\na \boldsymbol{\Phi}|^2$ and $(\lambda(\varphi_k)-\lambda(\varphi))^2|\delta|^2|\na \boldsymbol{\Phi}|^2$ gives
\begin{align*}
&\norm{(\eta(\varphi_k)-\eta(\varphi))\delta \na\boldsymbol{\Phi}}_{L^2(\mathbf{L}^2)}\ra 0\quad\text{ as }k\ra\infty,\\
&\norm{(\lambda(\varphi_k)-\lambda(\varphi))\delta \na\boldsymbol{\Phi}}_{L^2(\mathbf{L}^2)}\ra 0\quad\text{ as }k\ra\infty.
\end{align*}
Therefore, by the weak convergence $\bv_k\rightharpoonup \bv$ in $L^2(\mathbf{H}^1)$, $\div(\bv_k)\rightharpoonup \div(\bv)$ in $L^2(L^2)$ and $p_k\rightharpoonup p$ in $L^{\frac{4}{3}}(L^2)$, we easily deduce that
\begin{align}
\intT \intO \delta T(\bv_k,p_k)\colon\na\boldsymbol{\Phi}\d x\d t  &\ra \intT \intO\delta T(\bv,p)\colon\na\boldsymbol{\Phi}\d x\d t, \\
\intT\intO \delta\nu\bv_k\cdot\boldsymbol{\Phi}\d x\d t&\ra \intT\intO \delta\nu\bv\cdot\boldsymbol{\Phi}\d x\d t,
\end{align}
as $k\ra \infty$, where we used that $\delta\boldsymbol{\Phi}\in L^4(\mathbf{H}^1)$. Using $\varphi_k\ra\varphi$ strongly in $L^2(W^{1,3})$ and the continuous embedding $\mathbf{H}^1\hookrightarrow \mathbf{L}^6$, we have
\begin{equation*}
\intT\intO |\delta|^2|\boldsymbol{\Phi}|^2|\na\varphi_k-\na\varphi|^2\d x\d t 
\leq C\norm{\delta}_{L^{\infty}(0,T)}^2\norm{\boldsymbol{\Phi}}_{\mathbf{H}^1}^2\norm{\varphi_k-\varphi}_{L^2(W^{1,3})}^2\ra 0\quad\text{ as }k\ra\infty,
\end{equation*}
meaning $\delta\boldsymbol{\Phi}\cdot\na\varphi_k\ra \delta\boldsymbol{\Phi}\cdot\na\varphi$ strongly in $L^2(L^2)$. Again by the product of weak-strong convergence, it follows
\begin{equation}
\label{limit_23}\intT\intO \delta \mu_k\na\varphi_k\cdot\boldsymbol{\Phi}\d x\d t\ra \intT\intO \delta \mu\na\varphi\cdot\boldsymbol{\Phi}\d x\d t\quad\text{ as }k\ra\infty.
\end{equation}
By the specific form of $N_{\sigma,k}$ and since $\varphi_k\ra\varphi,~\sigma_k\ra\sigma$ strongly in $L^2(L^3)$, using a similar argument as for (\ref{limit_23}) yields
\begin{equation*}
N_{\sigma,k}\delta\boldsymbol{\Phi}\ra N_{\sigma}\delta\boldsymbol{\Phi}\quad\text{ strongly in }L^2(\mathbf{L}^2).
\end{equation*}
Consequently, by the product of weak-strong convergence we obtain
\begin{equation}
\label{limit_24}\intT\intO \delta N_{\sigma,k}\na\sigma_k\cdot\boldsymbol{\Phi}\d x\d t \ra \intT\intO \delta N_{\sigma}(\varphi,\sigma)\na\sigma\cdot\boldsymbol{\Phi}\d x\d t\quad\text{ as }k\ra\infty.
\end{equation}
Now we can pass to the limit in (\ref{limiting_equations_1})-(\ref{limiting_equations_5}) to obtain
\begin{align}
\nonumber \intT \delta(t)\langle\pa_t\varphi,w_j\rangle_{H^1,(H^1)^*}\d t &=\intT \delta(t)\left(\intO -m(\varphi)\na\mu\cdot\na w_j + \Gamma_{\varphi}w_j\d x\right)\d t\\
\label{limiting_equations_1a}&\quad  - \intT\delta(t)\left(\intO \na\varphi\cdot\bv w_j + \varphi\Gamma_{\bv}w_j \d x\right)\d t,\\
\label{limiting_equations_2a}\intT\intO \delta(t) \mu w_j\d x\d t &=\intT \intO \delta(t)(  \epsilon^{-1}\psi'(\varphi)w_j + \epsilon\na\varphi\cdot\na w_j - \chi_{\varphi}\sigma w_j) \d x\d t,\\
\nonumber \intT \delta(t)\langle \pa_t\sigma{,}w_j\rangle_{H^1,(H^1)^*}\d t&=\intT\intO \delta(t) (-n(\varphi)\na N_{\sigma}\cdot\na w_j -\Gamma_{\sigma}w_j)\d x\d t\\ 
\nonumber&\quad - \intT\intO \delta(t) (\na\sigma\cdot\bv w_j + \sigma\Gamma_{\bv}w_j)\d x\d t\\ 
\label{limiting_equations_3a}&\quad + \intT\delta(t)\left(\int_{\pa\Omega}b(\sigma_{\infty}-\sigma)w_j\d\mathcal{H}^{d-1}\right)\d t ,\\
\label{limiting_equations_4a}\intT \intO \delta(t) T(\bv,p)\colon \na\boldsymbol{\Phi}\d x\d t &= \intT\intO\delta(t) (-\nu\bv +\mu\na\varphi + N_{\sigma}\na\sigma)\cdot\boldsymbol{\Phi}\d x\d t,\\
\label{limiting_equations_5a}\intT\delta(t)\left( \intO \div(\bv)\Phi\d x\right)\d t &= \intT \delta(t)\left(\intO \Gamma_{\bv}\Phi\d x\right)\d t.
\end{align}
Since these equations hold for every $\delta\in C_0^{\infty}(0,T)$, we obtain that $\{\varphi,\sigma,\mu,\bv,p\}$ satisfies (\ref{weak_form_1}) with $\Phi = w_j$ for almost all $t\in (0,T)$ and all $j\geq 1$. Furthermore, (\ref{limiting_equations_5a}) implies $\div(\bv) = \Gamma_{\bv}(\varphi,\sigma)$ a.e. in $Q$. As $\{w_j\}_{j\in\N}$ is a basis for $H_N^2$ and $H_N^2$ is dense in $H^1$, this implies that $\{\varphi,\sigma,\mu,\bv,p\}$ satisfies (\ref{weak_form_1_eq_2})-(\ref{weak_form_1_eq_4}) for all $\Phi\in H^1$ and (\ref{weak_form_1_eq_1}) for all $\boldsymbol{\Phi}\in \mathbf{H}^1$.\newline\newline
\underline{\textbf{Step 6}:} We finally want to show that the initial conditions hold.
To this end, we notice that $\varphi_k(0)\ra\varphi_0$ in $L^2$. Furthermore, we know that $\varphi_k\ra\varphi$ strongly in $C^0([0,T];L^2)$, meaning $\varphi_k(0)\ra\varphi(0)$ strongly in $L^2$. But this already implies $\varphi(0)=\varphi_0$. Furthermore, since $\sigma$ belongs to $C^0([0,T];(H^1)^*)$, we see that $\sigma(0)$ is well-defined as an element of $(H^1)^*$. Furthermore, by the strong convergence $\sigma_k\ra\sigma$ in $C^0([0,T];(H^1)^*)$, we obtain for arbitrary $\zeta\in H^1$ that
\begin{equation*}
\lim_{n\ra \infty}\langle \sigma_k(0){,}\zeta\rangle_{H^1} = \langle\sigma(0){,}\zeta\rangle_{H^1}.
\end{equation*} 
By the weak convergence $\sigma_k(0)\ra\sigma_0$ in $L^2$, this implies
\begin{equation*}
\langle\sigma_0{,}\zeta\rangle_{H^1} = \lim_{n\ra \infty}\langle\sigma_k(0){,}\zeta\rangle_{H^1} = \langle\sigma(0){,}\zeta\rangle_{H^1},
\end{equation*}
which yields $\sigma(0) = \sigma_0$ in $(H^1)^*$. Finally, the energy inequality (\ref{estimate_solution_1}) follows from (\ref{central_estimate}) by weak (weak-star) lower-semicontinuity of the norms and dual norms. Having shown all these convergences, we proved the main result Theorem \ref{main_theorem}.

\section*{Acknowledgements} 
This work was supported by the RTG 2339 "Interfaces, Complex Structures, and Singular Limits" of the German Science Foundation (DFG). The support is gratefully acknowledged. Moreover, the authors are grateful to Prof. Helmut Abels for many interesting discussions
on topics related to this work.

\printbibliography[heading=bibintoc,title={References}]
\end{document}